\documentstyle{amsppt}
\mag=\magstep1
\pagewidth{30pc}
\pageheight{43pc}
\document
\topmatter
\title
Hilbert modular forms and $p$-adic Hodge theory
\endtitle
\author
Takeshi Saito
\endauthor
\affil
Department of Mathematical Sciences\\
University of Tokyo
\endaffil
\address
Tokyo 153-8914 Japan
\endaddress
\email
t-saito\@ ms.u-tokyo.ac.jp
\endemail
\endtopmatter

We consider the $p$-adic Galois representation associated
to a Hilbert modular form.
Carayol has shown that,
under a certain assumption,
its restriction to the local Galois group
at a place not dividing $p$
is compatible with the 
local Langlands correspondence [C2].
In this paper, 
we show that the same is true for
the  places dividing $p$,
in the sense of $p$-adic Hodge theory [Fo],
as is shown for an elliptic modular form in [Sa].
We also prove that the monodromy-weight
conjecture holds such representations.

We prove the compatibility
by comparing the $p$-adic and
$\ell$-adic representations
for it is already established
for $\ell$-adic representation [C2].
More precisely, we prove it
by comparing
the traces of Galois action and
proving the monodromy-weight conjecture.
The first task is to construct
the Galois representation
in purely geometric way in terms
of etale cohomology of an analogue of
Kuga-Sato variety and
algebraic correspondences
acting on it.
Then we apply the comparison theorem
of $p$-adic Hodge theory [Tj]
and weight spectral sequence [RZ], [M]
to compute the traces and monodromy operaters
in terms of the reduction modulo $p$.
We obtain the required equality
between traces by applying
Lefschetz trace formula
which has the same form for
$\ell$-adic and for cristalline
cohomology.
We deduce the monodromy-weight conjecture
from the Weil conjecture and
a certain vanishing of global sections.
The last vanishing result
is an analogue of the vanishing
of the fixed part $(\text{Sym}^{k-2}
T_\ell E)^{SL_2(\Bbb Z_\ell)}$
for $k>2$ for the universal elliptic curve $E$
over a modular curve in positive characteristic.

We briefly recall the
basic definitions on Hilbert modular forms
in Section 1 and
an $\ell$-adic representation
associated to it in Section 2.
The main compatibility result, Theorem 1,
and the monodromy-weight conjecture,
Theorem 2, are stated at the end of Section 2.
We recall a cohomological construction
of the $\ell$-adic representation 
in Section 3.
After introducing Shimura curves in Section 4
and recalling its modular interpretation
in Section 5,
we give a geometric construction 
of the $\ell$-adic representation 
in Section 6.
We extend the geometric construction
to semi-stable models in Section 7
and prove Theorems 1 and 2
in Section 8
admitting Proposition 1.
The last section 10 will be devoted to the proof
of Proposition 1.

The strategy of the proof is the same as in
the previous work in [Sa].
An essential part of the work consists of
understanding the work of Carayol [C1], [C2].
The author thanks Prof. K. Fujiwara for
the suggestion that the author's earlier proof for
totally real fields of odd degree should also work for
those of even degree
by using Carayol's construction of
$\ell$-adic representation.
He also thanks for Prof. F.Oort and Prof. A. de Jong
for teaching him sufficient conditions
for extension of abelian varieties.

Part of this work was done during
a stay at JAMI in Johns Hopkins University
in 1997,
a stay at IHP in $p$-adic semester
in 1997 and 
a stay at Paris-Nord in 1999.
The author would like to thank
their hospitality.

\medskip
\subheading{1. Hilbert modular form}

First, we briefly recall  
basic definitions on Hilbert modular forms
slightly modifying those in [Sh].
Let $F$ be a totally real number field of degree $g>1$
and $I=\{\sigma_1,\ldots,\sigma_g\}$
be the set of real embeddings $F\to \Bbb R$.
We fix a multiweight
$k=(k_1,\cdots,k_g,w)\in \Bbb N^{I+1}$ which is 
a $g+1$-uple of integers satisfying the conditions
$w\ge k_i\ge 2$ and $k_i\equiv w\bmod 2$.
The space $S^{(k)}_{\Bbb C}$
of cusp forms of multiweight $k$ is defined
as follows.

Let $X^I$
be the $g$-fold self product of the union 
$X=\Bbb P^1(\Bbb C)-\Bbb P^1(\Bbb R)$ of the upper 
and lower half planes.
It has a natural left action of
$GL_2(\Bbb R)^I$ by
$$\gamma(\tau)=\frac{a\tau+b}{c\tau+d}
=\left(\frac{a_i\tau_i+b_i}{c_i\tau_i+d_i}\right)_i\in X^I$$
for $\gamma=\pmatrix a&b\\c&d\endpmatrix=
\left(\pmatrix a_i&b_i\\c_i&d_i\endpmatrix\right)_i\in GL_2(\Bbb R)^I$
and $\tau=(\tau_i)_i\in X^I$.
It induces a left action of
$GL_2(F)$ on
$X^I\times GL_2(\Bbb A_{F,f})$
defined by $\gamma(\tau, g)=(\gamma(\tau),\gamma g)$.
Here $GL_2(F)$ is naturally embedded in
$GL_2(\Bbb A_{F})=
GL_2(\Bbb R)^I\times GL_2(\Bbb A_{F,f})$.
We also consider the right action
of $GL_2(\Bbb A_{F,f})$
on $X^I\times GL_2(\Bbb A_{F,f})$
defined by
$(\tau, g)g'=(\tau,gg')$.

A complex valued continuous function
$f=f(\tau,g)$ on $X^I\times 
GL_2(\Bbb A_{F,f})$ 
is said to be holomorphic if 
the function $\tau\mapsto f(\tau,g)$ on $X^I$
is holomorphic for each $g$ and
the map $g\mapsto (\tau\mapsto f(\tau,g))$
is locally constant on 
$GL_2(\Bbb A_{F,f})$.
The actions of $GL_2(F)$ and of
$GL_2(\Bbb A_{F,f})$ on holomorphic functions on
$X^I\times GL_2(\Bbb A_{F,f})$ are defined as follows.
For $\gamma\in GL_2(F)$ and a holomorphic function
$f$ on $X^I\times GL_2(\Bbb A_{F,f})$,
we define $\gamma^{(k)*}f =\gamma^*f$ to be 
$$
(\gamma^*f)(\tau,g)=
\frac{\det(\gamma)^{\frac{w+k-2}2}}{(c\tau+d)^k}f(\gamma \tau,\gamma g)
=
\prod_i\frac{\det(\gamma_i)^{\frac {w+k_i-2}2}}{(c_i\tau_i+d_i)^{k_i}}
 f(\gamma \tau,\gamma g).
$$
Here the image of $\gamma\in GL_2(F)$ in $GL_2(\Bbb R)^I$
is denoted by
$\left(\pmatrix a_i&b_i\\c_i&d_i\endpmatrix\right)_i$.
For $g'\in GL_2(\Bbb A_{F,f})$ and a holomorphic function
$f$ on $X^I\times GL_2(\Bbb A_{F,f})$,
we define $g^{\prime *}f$ by
$g^{\prime *}f(\tau, g)=f(\tau,gg')$.

For an open compact subgroup $K\subset GL_2(\Bbb A_{F,f})$,
a modular form of multiweight $k$ and of level $K$
is defined to be a holomorphic function $f$
on $X^I\times GL_2(\Bbb A_{F,f})$ invariant under
the action of $GL_2(F)$ and $K$.
Let 
$$M^{(k),K}_{\Bbb C}=
\left\{f\left|\matrix\text{ holomorphic function on 
$X^I\times GL_2(\Bbb A_{F,f})$
such that}\\
\text{ 
$\gamma^{(k)*}f=f$
for all $\gamma\in GL_2(F)$
and $g^*f=f$
for all $g\in K$}
\endmatrix\right.\right\}$$
be the space of
modular forms of multiweight $k$ and of level $K$.
We put
$M^{(k)}_{\Bbb C}=\bigcup_K M^{(k),K}_{\Bbb C}$
and call its element a modular form of multiweight $k$.

We recall the Fourier expansion of
a modular form and the definition of
the space of cusp forms.
Let $\psi_f:\Bbb A_{F,f}\to \Bbb C^\times$
be the finite part of the additive character
$$\Bbb A_{F}/F
@>{\text{Tr}_{F/\Bbb Q}}>>
\Bbb A/\Bbb Q
@<{\sim}<< (\hat {\Bbb Z}\times \Bbb R)/\Bbb Z
@>>>
\Bbb R/\Bbb Z
@>{a\mapsto\exp(2\pi\sqrt{-1} a)}>>\Bbb C^\times$$
and
$e_F$ be the function 
$$e_F(\tau, a)=\exp(2\pi\sqrt{-1}
\sum_i\tau_i)\cdot\psi_f (a)$$
on $X^I\times \Bbb A_{F,f}$.
Let $f$ be a modular form of multiweight $k$.
Then there exists a function
$c_z(\sigma,g,f)$ on $(z,\sigma,g)\in 
F \times \{\pm1\}^I\times GL_2(\Bbb A_{F,f})$
satisfying 
$$f\left(\tau, \pmatrix 1& b\\ 0&1\endpmatrix g\right)
=
\sum_{z\in F}
c_z(\text{sgn}(\text{Im }\tau),g,f)e_F(z\tau,zb)$$
since
$f\left(\tau, \pmatrix 1& b\\ 0&1\endpmatrix g\right)=
f\left(\tau+\beta, \pmatrix 1& b+\beta\\ 0&1 \endpmatrix g\right)$
for $\beta\in F$.
The Fourier coefficients
$c_z(\sigma,g,f)$ are 0
unless $z=0$ or $\text{sgn}(z)\sigma$
is totally positive.
In fact, 
a modular form is necessarily holomorphic at cusps
since we assume $[F:\Bbb Q]>1$.
A cusp form is defined to be a modular form $f$
satisfying $c_0(\sigma,g,f)=0$
for all $(\sigma,g)\in \{\pm1\}^I\times \Bbb A_{F,f}^\times$.
Unless $k_i$ are constant, 
we have $c_0(\sigma,g,f)=0$ and hence
a modular form is necessarily a cusp form.
In the following,
for a cusp form $f$,
we drop $\sigma$ in the notation and write
$c_z(g,f)=c_z(\text{sgn} z,g,f)$.
We say a cusp form $f$ is normalized
if $c_1(1,f)=1$.
Let $S^{(k)}_{\Bbb C}$
denote the space of
cusp forms of multiweight $k$.
For an open compact subgroup $K\subset GL_2(\Bbb A_{F,f})$,
we put
$S^{(k),K}_{\Bbb C}=
S^{(k)}_{\Bbb C}\cap M^{(k),K}_{\Bbb C}$.

We recall the definition of
the Dirichlet series $L(f,s)$ associated to
a cusp form $f$.
Let $D^{-1}=\{b\in F|\text{Tr}_{F/\Bbb Q}(O_Fb)\subset \Bbb Z\}$
be the codifferent ideal and
let $\hat T=\hat O_F\oplus D^{-1}\hat O_F\subset \Bbb A_{F,f}^2$
be a lattice.
For an integral ideal ${\frak n}\subset O_F$,
we define an open compact subgroup
$K_1({\frak n})\subset GL_{\hat O_F}(\hat T)
\subset GL_2(\Bbb A_{F,f})$ to be
$$\align
&K_1({\frak n})=
\left\{g\in GL_2(\Bbb A_{F,f})\left|g\hat T=\hat T, 
g\binom 10\equiv \binom10 \bmod\frak n\hat T
\right.\right\}\\
=&
\left\{\pmatrix a&b\\c&d\endpmatrix\in GL_2(\Bbb A_{F,f})
\left|
\matrix
a,d \in  \hat O_F, b\in D^{-1}\hat O_F,c\in D\hat O_F,
ad-bc\in \hat O_F^\times,\\
a\equiv 1\bmod \frak n\hat O_F,
c\equiv 0\bmod \frak n D\hat O_F
\endmatrix
\right.\right\}.
\endalign
$$
Let $f\in S^{(k),K_1(\frak n)}_{\Bbb C}$
be a cusp form.
For an idele
$d\in \Bbb A_{F,f}^\times$,
the Fourier coefficients
$c_1\left(\pmatrix 1& 0\\ 0&d\endpmatrix,f\right)$
depends only on the fractional ideal
$\frak m=(d^{-1})$ and
is 0 unless $\frak m\subset O_F$.
Here $(d^{-1})$ denotes the fractional ideal
$d^{-1}\hat O_F\cap F$.
We put 
$c(\frak m,f)=c_1\left(\pmatrix 1& 0\\ 0&d\endpmatrix,f\right)$
for an ideal $\frak m\subset O_F$
by taking an idele $d\in \Bbb A_{F,f}^\times$ such that
$\frak m=(d^{-1})$.
Since $c_z\left(\pmatrix1& 0\\0& d\endpmatrix,f\right)=
\prod_iz_i^{-\frac{w-k_i}2}
c_1\left(\pmatrix1& 0\\0& z^{-1}d\endpmatrix,f\right)$,
we have
$$f\left(\tau,\pmatrix 1&b\\0& d\endpmatrix\right)
=\sum_{\frak m\subset O_F}c(\frak m,f)
\sum_{z\in F^\times, 
zO_F=(d)\frak m,\text{sgn}z\tau=1}
\prod_iz_i^{\frac{w-k_i-2}2}e_F(z\tau,zd^{-1}b).$$
We define the Dirichlet series by
$$L(f,s)=\sum_{\frak m\subset O_F}c(\frak m,f)N\frak m^{1-s}.$$
It follows from the strong approximation theorem that
a cusp form $f\in S^{(k),K_1({ \frak n})}_{\Bbb C}$
is determined by the Fourier
coefficients $c(\frak m,f)$ i.e. by the L-series
$L(f,s)$.

We consider the natural action 
$g^*f(\tau,g')=f(\tau,g'g)$ of 
$g\in GL_2(\Bbb A_{F,f})$ on $S^{(k)}_{\Bbb C}$.
For $g\in GL_2(\Bbb A_{F,f})$
and open compact subgroups $K,K'\subset GL_2(\Bbb A_{F,f})$
satisfying
$g^{-1}Kg\subset K'$,
the map $g^*$ sends $S^{(k),K'}_{\Bbb C}$ into
$S^{(k),K}_{\Bbb C}$ :
$g^*:S^{(k),K}_{\Bbb C}\to
S^{(k),K'}_{\Bbb C}$.
For open compact subgroups $K'\subset K\subset GL_2(\Bbb A_{F,f})$,
we have the trace map
$\text{Tr}:S^{(k),K'}_{\Bbb C}\to
S^{(k),K}_{\Bbb C}$
defined by
$f\mapsto \sum_{g\in K/K'}g^*f$.
Hence by taking the fixed part,
we recover
$S^{(k),K}_{\Bbb C}=(S^{(k)}_{\Bbb C})^K$.
For an open compact subset 
$T\subset GL_2(\Bbb A_{F,f})$
stable under the actions of
an open compact subgroup $K\subset GL_2(\Bbb A_{F,f})$
on the bothside,
the action
$T^*:
S^{(k),K}_{\Bbb C}\to
S^{(k),K}_{\Bbb C}$
is defined by
$f\mapsto \sum_{g\in T/K}g^*f$.
If $T=KgK$ for $g\in GL_2(\Bbb A_{F,f})$,
it is the same as the composite
$$S^{(k),K}_{\Bbb C}@>{g^*}>>
S^{(k),K\cap gKg^{-1}}_{\Bbb C}@>{\text{Tr}}>>
S^{(k),K}_{\Bbb C}.$$

We define the Hecke operators.
Let $\frak p$ be a maximal ideal of $O_F$.
First we consider the case where
an open compact subgroup 
$K\subset GL_2(\Bbb A_{F,f})$
is the product $K=K_{\frak p}\times K^{\frak p}$
of $K_{\frak p}=GL_{O_{F_{\frak p}}}(T)$
 for some $O_{F_{\frak p}}$-lattice $T$ of
$F_{\frak p}^2$ with a prime-to-$\frak p$
part $K^{\frak p}$.
We put 
$$\align
T_{\frak p}=&\
\{g\in GL_2(F_{\frak p})|
gT\supset T, gT/T\simeq O_F/\frak p\}\times K^{\frak p},
\\
R_{\frak p}=&\ 
\{g\in GL_2(F_{\frak p})|
gT\supset T, gT/T\simeq (O_F/{\frak p})^{\oplus 2}\}
\times K^{\frak p}.
\endalign$$
They define endomorphisms
also denoted by $T_{\frak p},R_{\frak p}$ on 
$S^{(k),K}_{\Bbb C}$ respectively.
They are independent of the choice of
$K_{\frak p}$-stable lattice $T$
and depends only on $K$.
If $K$ is of the form $K_{\frak q}K^{\frak q}$ as above for 
another prime ${\frak q}$,
the operaters $T_{\frak p},R_{\frak p},
T_{\frak q},R_{\frak q}$ are commutative to each other.
Note that $T_{\frak p},R_{\frak p}$ defined above
are closely related to but slightly different from 
$T(\frak p),T(\frak p,\frak p)$ in [Sh].
It is analogous to those in [D1].

When $K=K_1(\frak n)$,
for an arbitrary integral ideal $\frak m\subset O_F$,
we put 
$$\align
T_{\frak m}=&
\{g\in GL_2(\Bbb A_{F,f})
|g(\hat O_F\oplus D\hat O_F)\supset 
\hat O_F\oplus D\hat O_F, \det g^{-1}\hat O_F=
\frak m\hat O_F\}\\
R_{\frak m}=&\cases
\left\{\left.
\pmatrix a^{-1}&0\\0&a^{-1}
\endpmatrix \right|a\in \Bbb A_{F,f}^\times,\ a\hat O_F=
\frak m\hat O_F\right\}
\cdot K_1(\frak n)
&\quad\text{ if $(\frak n,\frak m)=1$,}\\
\emptyset 
&\quad\text{ if otherwise.}
\endcases
\endalign$$
and let $T_{\frak m}, R_{\frak m}$ also denote the endomorphism
on $S^{(k),K_1(\frak n)}_\Bbb C$.
When $\frak m=\frak p\nmid \frak n$, 
the two definitions are the same.
The Hecke operaters satisfy the formal equality
$$\sum T_{\frak m}N\frak m^{-s}=
\prod_{\frak p}
\left(1-T_{\frak p}N{\frak p}^{-s}+
R_{\frak p}N{\frak p}^{1-2s}\right)^{-1}.$$
For a cusp form $f\in S_{\Bbb C}^{(k),K_1(\frak n)}$,
by an elementary computation, we see
$c(1,T_{\frak m}(f))=c({\frak m},f)\cdot N\frak m $.
Suppose a cusp form $f\in S_{\Bbb C}^{(k),K_1(\frak n)}$
is an eigenform for all the Hecke operators
$T_{\frak m},R_{\frak m}$ and 
let $\chi_f(\frak m)$ be the eigenvalue of
$R_{\frak m}$.
Then since $f=0$ if $c(\frak m,f)=0$ for all $\frak m$,
we see $c(1,f)\neq 0$ and the eigenvalue of
$T_{\frak m}$ is $c(\frak m,f)\cdot N\frak m /c(1,f)$.
Hence if $f$ is further normalized,
the $L$-series has the Euler product
$$L(f,s)=
\prod_{\frak p}
\left(1-c(\frak p,f)N\frak p^{1-s}+
\chi_f(\frak p)N\frak p^{1-2s}\right)^{-1}.$$
The local Euler factor is
$L_{\frak p}(f,T)=
1-c(\frak p,f)N\frak p\ T+\chi_f(\frak p)N\frak p\ T^2.$

An irreducible subrepresentation
$\pi$ of the  representation
$S^{(k)}_{\Bbb C}$ of the adele group 
$GL_2(\Bbb A_{F,f})$
is called a cuspidal automorphic representation
of multiweight $k$.
It is known that we have a direct sum decomposition
$S^{(k)}_{\Bbb C}=\bigoplus_{\pi}\pi$
where $\pi$ runs cuspidal automorphic representations
of multiweight $k$.
For a cuspidal automorphic representation $\pi$,
there exists a largest ideal
$\frak n$, called the level of $\pi$,
such that the $K_1(\frak n)$-fixed part
$\pi^{K_1(\frak n)}=\pi\cap S^{(k),K_1(\frak n)}_{\Bbb C}$
is non-zero.
If $\frak n$ is the level of a cuspidal automorphic representation
$\pi$,
the dimension of the space
$\pi\cap S^{(k),K_1(\frak n)}_{\Bbb C}$ is 1.
A non-zero element in this space is
an eigenform for all the Hecke operators
$T_{\frak m}, R_{\frak m}$.
Hence, it is generated by a normalized cusp form $f$.
The Fourier coefficients $c(\frak m,f)$ is determined 
by the condition that $c(\frak m,f)N\frak m$ is equal to
the eigenvalue of the Hecke operater $T_{\frak m}$. 
We call such $f$ a normalized eigen new form.
Since the irreducible representation $\pi$ is
generated by $f$,
the correspondence $\pi\leftrightarrow f$ 
between the cuspidal automorphic representations
and the normalized new eigenforms is one-to-one.
In particular,
each irreducible factor $\pi$ appears
only once in the direct sum decomposition above.
In the following,
we let $\pi_f$
denote the cuspidal automorphic representation
generated by $f$ for a normalized new eigenform $f$.

Since a cuspidal automorphic representation
$\pi_f$ is irreducible,
the center 
$\Bbb A_{F,f}^\times\subset GL_2(\Bbb A_{F,f})$
acts on $\pi_f$ by the so-called central character
$\chi_{\pi_f}:\Bbb A_{F,f}^\times\to \Bbb C^\times.$
The conductor of the central character $\chi_{\pi}$
divides the level of $\pi_f$.
For $\frak p\nmid \frak n$,
we have $\chi_f(\frak p)=\chi_{\pi_f}(\frak p)^{-1}$.
Since $\chi_f|_{F^\times}=
\chi_{\pi_f}^{-1}|_{F^\times}=
N_{F/\Bbb Q}^{w-2}:F^\times\to \Bbb C^\times$,
the character 
$\chi_f$ is an algebraic Hecke character 
whose algebraic part is $N_{F/\Bbb Q}^{w-2}$.
Hence there is a character
$\epsilon_f:\Bbb A_{F,f}^\times/F^\times\to \Bbb C^\times$
of finite order and of conductor dividing 
$\frak n$
such that
$\chi_f=N_{F/\Bbb Q}^{w-2}\cdot \epsilon_f$.
Therefore, for $\frak p\nmid \frak n$,
the Euler factor is given by
$$
L_{\frak p}(f,T)=
1- c(\frak p,f)\cdot N \frak p\ T+
\epsilon_f(\frak p)\cdot N\frak p^{w-1}\ T^2.$$
Note that it is slightly different from the
Euler factor $L_{\frak p}(\pi_f,T)$.
We take this definition in order to make
the formula (2.1) simpler.

To define an $L$-structure 
$S_L^{(k)}$ of $S_{\Bbb C}^{(k)}$
over an number field $L$,
we recall the description of
the space $S^{(k),K}_{\Bbb C}$ of cusp forms
in terms of automorphic bundles [Mi] Chap.III.
Let $S=(S_K)_K$ be the canonical model of the Hilbert modular variety.
It is the canonical model $Sh(G,X)$
over the reflex field $\Bbb Q$
of the Shimura variety defined 
by $G=GL_{2,F}$ regarded as
an algebraic group over $\Bbb Q$
and the $G(\Bbb R)$-conjugacy class $X^I$ of
the homomorphism
$$\matrix
h:&\Bbb C^\times&\to
&G(\Bbb R)=
GL_2(\Bbb R)^I,\\
&a+b\sqrt{-1}&\mapsto &
\left(\pmatrix a&b\\-b&a\endpmatrix,\ldots,
\pmatrix a&b\\-b&a\endpmatrix\right)
\endmatrix$$
(cf.\ loc.\ cit.).
For an open compact subgroup $K\subset GL_2(\Bbb A_{F,f})$,
the set of complex points $S_K(\Bbb C)$ is given as the double cosets
$$S_K(\Bbb C)=
GL_2(F) \backslash X^I\times GL_2(\Bbb A_{F,f})/K.$$
The projective system
$S=(S_K)_K$
has an natural action of the group $GL_2(\Bbb A_{F,f})$.

We define an invertible sheaf $\omega^{(k)}$
on $S$ as an automorphic vector bundle $\Cal V(\Cal J^{(k)})$ as
follows.
Let $\check X^I=(\Bbb P^1_{\Bbb C})^I$ be the compact dual of $X^I$.
It has a natural action of $G_{\Bbb C}=(GL_{2,\Bbb C})^I$.
Let $\omega$ be the dual of the tautological quotient bundle
on $\Bbb P^1_{\Bbb C}$.
The line bundle
$\omega$ has a natural equivariant $GL_2$-action.
We define an equivariant $G_{\Bbb C}$-action on
the line bundle
$$\Cal J^{(k)}=
\bigotimes_{i\in I} pr_i^*\omega^{\otimes k_i-2}$$
on $\check X^I$ as follows.
For each $i\in I$, we define an equivariant action of
$GL_{2,\Bbb C}$ on $\omega^{\otimes k_i-2}$
to be the $\det^{-\frac {w-k_i}2}$-times
the tensor product of the natural action $\omega$.
By taking the tensor product,
we define a $G_{\Bbb C}$-equivariant action on $\Cal J=\Cal J^{(k)}$.
Let $G^c$ denotes the quotient of $G$ by
$\text{Ker}(N_{K/F}:F^\times \to \Bbb Q^\times)$.
Here $F^\times$ denotes the center of $GL_{2,F}$ as an algebraic group over
$\Bbb Q$.
Since the center $(\Bbb G_{m,\Bbb C})^I\subset G_{\Bbb C}$
acts by the $-(w-2)$-nd power
of the product character,
it defines a $G^c_{\Bbb C}$-equivariant bundle.
Hence as in [Mi] Chap.III,
we obtain a $G(\Bbb A_f)$-equivariant vector bundle
$\omega^{(k)}=\Cal V(\Cal J^{(k)})$ on 
$S(\Bbb C)=(S_K(\Bbb C))_K$.

We say an open compact subgroup $K\subset GL_2(\Bbb A_{F,f})$
is sufficiently small
if the following two conditions are satisfied.

(1)
The quotient $(gKg^{-1}\cap GL_2(F))/(gKg^{-1}\cap F^\times)$
does not have non-trivial element of
finite order
for all $g\in GL_2(\Bbb A_{F,f}).$

(2)
$N_{F/\Bbb Q}(K\cap F^\times)^{w-2}=1$.

\noindent
If $K$ satisfies the condition (1),
the canonical map $S_{K'}\to S_K$ is etale
for an open subgroup $K'\subset K$
and $S_K$ is smooth over $\Bbb Q$.
If $K$ satisfies the condition (2),
the invertible sheaf $\omega^{(k)}$ is defined on
$S_K(\Bbb C)$.
Following the definition, 
it is straightforward to check that
the space $M_{K,\Bbb C}^{(k)}$ 
of modular forms is identified with the
space of global sections 
$\Gamma(S_K(\Bbb C),\omega_{S_K}\otimes \omega^{(k)})$
for sufficiently small $K\subset GL_2(\Bbb A_{F,f})$.
Here $\omega_{S_K}$ denotes the canonical invertible sheaf
$\Omega^g_{S_K}$.

Let $L\subset \Bbb C$ be
a number field which contains all the conjugates of $F$.
Then the $G^c$-equivariant
bundle $\Cal J^{(k)}$ is defined over $L$.
The $G(\Bbb A_f)$-equivariant invertible sheaf $\omega^{(k)}$
has the canonical model $\omega^{(k)}_L$
defined over $L$  by [Mi] Chap.III Theorem 5.1.(a).
Hence $M_K^{(k)}(\Bbb C)$
has a natural $L$-structure
$M_K^{(k)}(L)=
\Gamma(S_K,\omega_{S_K}\otimes \omega^{(k)}_L)$.
Further we put
$S_K^{(k)}(L)=
M_K^{(k)}(L)\cap S_K^{(k)}(\Bbb C)$
and get
$S_K^{(k)}(\Bbb C)=S_K^{(k)}(L)\otimes_L\Bbb C$.
The last equality is deduced from the fact that the
Fourier expansion is defined algebraically using
the HBAV-analogue of Tate curve as in [Ka].

We describe the invertible sheaf
$\omega^{(k)}$ in terms of the moduli problem.
The Hilbert modular variety $S_K$ is the coarse moduli
scheme for the following functor:
To a scheme $S$ over $\Bbb Q$,
associate the set of isomorphism classes
of 
abelian schemes $A$ over $S$ upto isogeny of dimension g
endowed with a ring homomorphism
$F\to \text{End}_S(A)$
such that $\text{Lie}\ A$ is an invertible
$F\otimes_{\Bbb Q} O_S$-module,
together with a weak polarization respect to $F$ and
with a level structure modulo $K$.
Although it is not representable,
for $K$ sufficiently small,
there is an abelian scheme $a:A\to S_K$ upto isogeny
with multiplication by $F$,
which is almost the universal abelian scheme.
The cotangent bundle $\omega=a_*\Omega^1_{A/S_K}$
is an invertible $O_{S_K}\otimes F$-module
and its dual is $\text{Lie}\ A$.
The higher direct image of the relative de Rham complex
$H=R^1a_* \Omega^\bullet_{A/S_K}$
is a locally free
$O_{S_K}\otimes_{\Bbb Q}F$-module of rank 2.
By $F\otimes_{\Bbb Q}L=\prod_{i\in I}L$,
we have direct sum decompositions
$\omega\otimes_{\Bbb Q} L
=\bigoplus_i\omega_i$ and 
$H\otimes_{\Bbb Q} L=\bigoplus_iH_i.$
Here
$\omega_i=\omega\otimes_{F,\sigma_i}L,
H_i=H\otimes_{F,\sigma_i}L$
are locally free $O_{S_K}\otimes_{\Bbb Q}L$-modules
of rank 1 and 2 respectively.
Then we have
$$\omega^{(k)}_L=
\bigotimes_i\left((\Lambda^2 H_i)^{\otimes \frac{w-k_i}2}\otimes
\omega_i^{\otimes k_i-2}\right).$$

The representation $S^{(k)}_{\Bbb C}$ of
$GL_2(\Bbb A_{F,f})$ has an
$L$-structure $S^{(k)}_L=\varinjlim_K S^{(k),K}_L$.
For a normalized eigen new form $f$,
we see that the field
$L(f)=L(c({\frak m},f),{\frak m}\subset O_F)$ generated by the
Fourier coefficients
is of finite degree over $L$
and that the representation $\pi_f$
has an $L(f)$-structure
$\pi_{f,L(f)}$.
In fact, let ${\frak n}$ be the level of $f$
and consider the Hecke algebra $T({\frak n})_L=
L[T_{\frak m},{\frak m}\subset O_F]
\subset \text{End}_L(S^{(k),K_1({\frak n})}_L)$.
Then $L(f)$ is the image of $T({\frak n})\to \Bbb C$
defined by the action on the subspace of
$S_{\Bbb C}^{(k),K_1({\frak n})}$ by $f$
and is of finite degree.
The intersection
$S_{L(f)}^{(k),K_1({\frak n})}\cap \pi_f$ is
identified with
$\text{Hom}_{T({\frak n})_{L(f)}}
(L(f),S_{L(f)}^{(k),K_1({\frak n})})$
and is of dimension one over $L(f)$.
The subrepresentation of $S^{(k)}_{L(f)}$
generated by this line gives
an $L(f)$-structure
$\pi_{f,L(f)}$.
We have direct sum decomposition
$S^{(k)}_L=\bigoplus_f\pi_{f,L(f)}$
where $f$ runs the conjugacy classes of
eigen newforms $f$
over $L$.
The Euler factors
$L_{\frak p}(f,T)$ have the coefficients in
the number field $L(f)$.

\medskip
\subheading{2. $\ell$-adic representation
associated to a Hilbert modular form: Main results}

We recall the definition of the $\ell$-adic representation
associated to a Hilbert modular form.
Let $F$ be a totally real number field and
let $f$ be a normalized new eigen form of multiweight $k$.
Let $L\subset \Bbb C$ be 
a number field which contains all the conjugate
of $F$
and take a finite place $\lambda$ of the number field $L(f)$
generated by the Fourier coefficients $c({\frak m},f)$.
Then an $\ell$-adic representation
$\rho:G_F\to GL_2(L(f)_\lambda)$
is said to be associated to $f$,
if the following condition is satisfied 
for almost all finite places ${\frak p}$ of $F$.

The representation $\rho$ is unramified at ${\frak p}$
and the eigen polynomial of the geometric Frobenius
$Fr_{\frak p}$ is equal to the Euler factor $L_{\frak p}(f,T)$
$$\det (1-\rho(Fr_{\frak p})T)=L_{\frak p}(f,T)
=1-c({\frak p},f)\cdot N{\frak p} \ T+
\epsilon_f({\frak p})\cdot N{\frak p}^{w-1}\ T^2.
\tag 2.1$$

\noindent
In practice, the finite subset of places to be omitted
consists of those dividing the level ${\frak n}$ of $f$
or the prime $\ell$ below $\lambda$.
The existence is established by an accumulation
of works of many people
[O], [C2], [RT],[BR], [Ta1].
Since it is known to be irreducible
[Ta2] Proposition 3.1,
Chebotarev density implies the uniqueness.
In the following, we recall a theorem of Carayol [C2]
which asserts not only the existence but
also gives a precise description of the restriction to
the decomposition group
$D_{\frak p}=\text{Gal}(\bar F_{\frak p}/F_{\frak p})$
at finite places $p\nmid\ell$ including those dividing
the level $\frak n$.
The description is given
in terms of local Langlands correspondence,
recalled in the following.

Let $\pi_{f,L(f)}$ be the $L(f)$-structure
of the cuspidal automorphic representation
of $GL_2(\Bbb A_{F,f})$ associated to $f$.
Let $\pi_{f,L(f)}=\bigotimes_{\frak p}\pi_{f,L(f),{\frak p}}$
be the factorization 
into the tensor product
of irreducible admissible representations $\pi_{f,L(f),{\frak p}}=
\pi_{f,L(f)}^{K_1(n)^{\frak p}}$
of $GL_2(F_{\frak p})$ over $L(f)$.
Here ${\frak n}$ is the level of $f$
and $K_1({\frak n})^{\frak p}$ is the prime-to-${\frak p}$
component of $K_1({\frak n})=
K_1({\frak n})_{\frak p}\cdot K_1({\frak n})^{\frak p}$.
To attach an $L(f)$-rational
representation of Weil-Deligne group to
the $L(f)$-representation 
$\pi_{f,L(f),{\frak p}}$ of $GL_2(F_{\frak p})$,
we briefly recall the local Langlands correspondence.

To an irreducible admissible
representation $\pi$ of $GL_2(F_{\frak p})$,
the local Langlands correspondence
associates an $F$-semi-simple 
representation $\sigma(\pi)$ of the Weil-
Deligne
group ${}'W(\bar F_{\frak p}/F_{\frak p})$ of degree 2.
An $F$-semi-simple 
representation of the Weil-Deligne
group is a pair of
a semi-simple representation $(\rho,V)$
of the Weil group $W(\bar F_{\frak p}/F_{\frak p})$
with open kernel and a nilpotent endomorphism
$N$ of $V$
satisfying
$\rho(\sigma)N\rho(\sigma)^{-1}=N{\frak p}^{n(\sigma)}N$.
Here $N{\frak p}$ is the norm of ${\frak p}$
and $n:W(\bar F_{\frak p}/F_{\frak p})\to \Bbb Z$
is the canonical surjection sending
a geometric Frobenius in
$W(\bar F_{\frak p}/F_{\frak p})$ to 1.
A representation $(\rho',N)$ of the Weil-Deligne group
is called unramified if $\rho'$ is unramified
and $N=0$.

Among several ways to normalize the local Langlands
correspondence,
here we consider the so-called Hecke correspondence
: $ \pi\mapsto \sigma_h(\pi)$ [De].
To describe the normalization,
we give the definition
for a spherial representation $\pi$ of $GL_2(F_{\frak p})$.
Let $K=GL_2(O_{\frak p})$
and let $\pi_{\frak p}$ be a prime element of $F_{\frak p}$.
We define
$\tau_{\frak p}$ and $\rho_{\frak p}$ to be the eigenvalues of
the action of the double cosets
$$\align
T_{\frak p}=&
\{g\in GL_2(F_{\frak p})|gO_{F_{\frak p}}^2\supset 
O_{F_{\frak p}}^2,
\det g O_{F_{\frak p}}={\frak p}^{-1}\}=
K\pmatrix \pi_{\frak p}^{-1}&0\\0&1\endpmatrix K,\\
R_{\frak p}=&
\{g\in F_{\frak p}^\times|
g O_{F_{\frak p}}={\frak p}^{-1}\}\cdot GL_2(O_{F_{\frak p}})=
K\pmatrix \pi_{\frak p}^{-1}&0\\0&\pi_{\frak p}^{-1}\endpmatrix K
\endalign$$
respectively on the 1-dimensional space
$\pi^K$.
Then the representation $\sigma_h(\pi)$
is the unramified semi-simple representation
characterized by
$$\det(1-Fr_{\frak p}T:\sigma_h(\pi))=
1-N{\frak p}^{-1}\rho_{\frak p}^{-1}\tau_{\frak p}\ T
+N{\frak p}^{-1}\rho_{\frak p}^{-1}\ T^2.$$
It is the same as to require, for the dual representation
$\check\sigma_h(\pi)$,
$$\det(1-Fr_{\frak p}T:\check\sigma_h(\pi))=
1-\tau_{\frak p}\ T+\rho_{\frak p}\cdot N{\frak p}\ T^2.$$
If $\pi$ is defined over a field $L$ of characteristic 0,
the representation $\sigma_h(\pi)$ is
defined over an algebraic closure $\bar L$
and the isomorphism class is invariant under
the Galois group $\text{Gal}(\bar L/L)$.
In other word, the representation
$\sigma_h(\pi)$ is $L$-rational,
not necessarily realized over $L$.

We apply the construction $\pi\mapsto \sigma_h(\pi)$ to
the local component
$\pi_{f,{\frak p}}$ of a cuspidal automorphic representation.
Thus, we obtain
an $F$-semi-simple $L(f)$-rational representation
$\check\sigma_h(\pi_{f,{\frak p}})$ of the
Weil-Deligne group ${}'W(\bar F_{\frak p}/F_{\frak p})$.
For a prime ${\frak p}\nmid {\frak n}(f)$,
it is an unramified representation.
The equality (1) above is then rephrased as
$$\det(1-\rho(Fr_{\frak p})T)=
\det(1-Fr_{\frak p}T:\check\sigma_h(\pi_{f,{\frak p}})).$$
It is the same as to say that
the semi-simplification of the unramified
representation
$\rho_{f,\lambda}|_{W_{\frak p}}$
is isomorphic to 
$\check\sigma_h(\pi)$.

On the other hand, 
to an $\ell$-adic representation 
of the local Galois group
$G_{\frak p}=\text{Gal}(\bar F_{\frak p}/F_{\frak p})$,
we attach a representation 
of the Weil-Deligne group 
${}'W(\bar F_{\frak p}/F_{\frak p})$.
First we consider the case where
${\frak p}\nmid \ell$.
Let $L_\lambda$ be a finite extension of
$\Bbb Q_\ell$.
Let $\rho:G_{\frak p}\to GL_{L_\lambda}(V)$
be a continuous $\ell$-adic representation.
Take a lifting $F\in W(\bar F_{\frak p}/F_{\frak p})$ of 
the geometric Frobenius and 
an isomorphism $\Bbb Z_{\ell}(1)\to\Bbb Z_{\ell}$
and identify them.
Let $t_\ell:I_{\frak p}\to\Bbb Z_{\ell}(1)\to\Bbb Z_{\ell}$
be the canonical surjection.
Then, by the monodromy theorem of Grothendieck,
there is a representation
${}'\rho=(\rho',N)$ of the Weil-Deligne group 
${}'W(\bar F_{\frak p}/F_{\frak p})$
characterized by the condition
$$\rho(F^n\sigma)=\rho'(F^n\sigma)\exp(t_\ell(\sigma) N)$$
for $n\in \Bbb Z$ and $\sigma\in I_{\frak p}$.
The isomorphism class of the representation $(\rho',N)$
of the Weil-Deligne group is independent of the
choice of the lifting $F$ or the isomorphism
$\Bbb Z_{\ell}(1)\to\Bbb Z_{\ell}$
and is determined by $\rho$.

For an $\ell$-adic representation $\rho$ of
$\text{Gal}(\bar F/F)$,
let $\rho_{\frak p}$
denote the restriction to
$\text{Gal}(\bar F_{\frak p}/F_{\frak p})$.
Let ${}'\rho_{\frak p}$
denote the representation of
the Weil-Deligne group attached to
$\rho_{\frak p}$ and 
let
${}'\rho_{\frak p}^{\text{F-ss}}$
denote its $F$-semi-simplification.

\proclaim{Theorem 0}{\rm [C2]}
Let $f$ be a normalized eigen newform of multiweight $k$
and
$\lambda|\ell$ be a finite place of
the number field $L(f)$.
We assume that,
if the degree $g=[F:\Bbb Q]$ is even,
there exists a finite place ${\frak v}$
such that the ${\frak v}$-factor 
$\pi_{f,{\frak v}}$ lies in the discrete series.
Then there exists an $\ell$-adic representation
$$\rho=\rho_{f,\lambda}:
G_F@>>>GL_{L(f)_\lambda}(V_{f,\lambda})$$
satisfying the following property:

\noindent
For a finite place ${\frak p}\nmid \ell$,
there is an isomorphism
$${}'\rho_{f,\lambda,{\frak p}}^{\text{F-ss}}\simeq
\check\sigma_h(\pi_{f,{\frak p}})$$
of representations of
the Weil-Deligne group
${}'W(\bar F_{\frak p}/F_{\frak p})$.
\endproclaim

\demo{Remark}
Since the right hand side is $L(f)$-rational,
Theorem implies that so is the left hand side.
For ${\frak p}\nmid {\frak n}(f)\ell$,
the isomorphism means that we have an equality
$$\det(1-Fr_{\frak p}T:\rho_{f,\lambda})=
\det(1-Fr_{\frak p}T:\check\sigma_h(\pi))
=L_{\frak p}(f,T).$$
Hence $V_{f,\lambda}$ in Theorem 1 is
the $\ell$-adic representation 
associated to $f$.
\enddemo

In this paper, we study the case ${\frak p}$ divides $\ell$.
Let $p$ be the characteristic of
a finite place ${\frak p}$ of $F$.
Let $F_{\frak p,0}$ denoted the maximal
unramified subfield in $F_{\frak p}$ .

We describe the construction attaching a representation of
Weil-Deligne group to a $p$-adic representation
of the local Galois group due to Fontaine [Fo].
Let $B_{st}$ be the ring defined by Fontaine.
It is an $\widehat{F_{\frak p,0}^{nr}}$-algebra
and admits a natural action of
the absolute Galois group $G_{F_{\frak p}}$,
a semi-linear action of the Frobenius
$\varphi$ and an action of the monodromy operator $N$.
For an open subgroup $J\subset I$
of the inertia,
the fixed part
$B_{st}^J$
is the completion
$\widehat{F_{\frak p,0}^{nr}}$
of a maximal unramified extension 
of $F_{\frak p,0}$.
In this paper, we neglect the filtration.
Let $L_\mu$ be a finite extension of $\Bbb Q_p$
and consider a continuous $p$-adic representation
$\text{Gal}(\bar F_{\frak p}/F_{\frak p})\to 
GL_{L_\mu}(V)$ of finite degree.
Let $\widehat{L_\mu^{nr}}$ denote the completion of
the maximum unramified extension of $L_\mu$.
We choose an arbitrary factor of 
$\widehat{F_{\frak p,0}^{nr}}
\otimes_{\Bbb Q_p} L_{\mu}$.
It is the same thing as to
fix an embedding
$\widehat{F_{\frak p,0}^{nr}}
\to \widehat{L_\mu^{nr}}$.
For an $L_\mu$-representation
$G_{F_{\frak p}}\to GL_{L_\mu}(V)$
of finite degree,
we put
$$D(V)=D_{pst}(V)
=\bigcup_{J\subset I}
(B_{st}\otimes V)^J
\otimes_{(\widehat{F_{\frak p,0}^{nr}}\otimes_{\Bbb Q_p}L_\mu)}
\widehat{L_\mu^{nr}}.$$
Here $J$ runs the open subgroups of the
inertia subgroup $I=I_{\frak p}$
and ${}^J$ denotes the $J$-fixed part.
The union
$\bigcup_{J\subset I}
(B_{st}\otimes V)^J$
is a $\widehat{F_{\frak p,0}^{nr}}
\otimes_{\Bbb Q_p}L_\mu$-module
since $B_{st}^J=\widehat{F_{\frak p,0}^{nr}}$.
It is known that $D(V)$ is
an $\widehat{L_\mu^{nr}}$-vector space of finite dimension
and
$\dim_{\widehat{L_\mu^{nr}}} D(V)\le \dim_{L_\mu} V$.
We say $V$ is potentially semi-stable (pst for short)
if we have an equality
$\dim_{\widehat{L_\mu^{nr}}} D(V)=\dim_{L_\mu} V$.

For a pst-representation $V$,
Fontaine defines a natural representation [Fo]
on $D(V)$ of
the Weil-Deligne group
${}'W(\bar F_{\frak p}/F_{\frak p})$
as follows [Fo].
By the Galois actions on $B_{st}$ and on $V$,
the quotient $G_F/J$ acts on
the $J$-fixed part
$(B_{st}\otimes V)^J$
for normal $J\subset G_F$.
Passing to the limit,
we obtain an action of
$G_F$ acts on 
the 
$\widehat{F_{\frak p,0}^{nr}}\otimes_{\Bbb Q_p}L_\mu$-module
$\bigcup_{J\subset I}
(B_{st}\otimes V)^J$.
The kernel is open in the inertia $I_{\frak p}$.
This Galois action is semi-linear with respect to
its natural action on
$\widehat{F_{\frak p,0}^{nr}}$
and the trivial action on $L_\mu$.
We modify it by using the Frobenius
$\varphi$ to get a
$\widehat{F_{\frak p,0}^{nr}}\otimes_{\Bbb Q_p}L_\mu$-linear
action
of the Weil group
$W(\bar F_{\frak p}/F_{\frak p})$
as follows.

Let $\Bbb F_{\frak p}$ denote the residue field of $\frak p$.
Recall that the Weil group
$W(\bar F_{\frak p}/F_{\frak p})$ is the inverse image
of the inclusion $\Bbb Z \to 
\text{Gal}(\overline {\Bbb F_{\frak p}}/
\Bbb F_{\frak p})$
sending 1 to the geometric Frobenius $Fr_{\frak p}$
by the canonical map
$\text{Gal}(\bar F_{\frak p}/F_{\frak p})\to
\text{Gal}(\overline {\Bbb F_{\frak p}}/\Bbb F_{\frak p})$.
Let $n:W(\bar F_{\frak p}/F_{\frak p})\to\Bbb Z$
be the canonical map and $q=p^f=N\frak p$.
Then by let $\sigma\in
W(\bar F_{\frak p}/F_{\frak p})$
act on $D(V)$
by $(\varphi^{f\cdot n(\sigma)}\otimes 1)\circ
\sigma\otimes \sigma$,
we get a
$\widehat{F_{\frak p,0}^{nr}}\otimes_{\Bbb Q_p}L_\mu$-linear action.
Taking the $\widehat{L_\mu^{nr}}$-component,
we obtain an
$\widehat{L_\mu^{nr}}$-linear representation $D(V)$
of the Weil group
$W(\bar F_{\frak p}/F_{\frak p})$.
The monodromy operator $N$
on $B_{st}$ induces
an $\widehat{L_\mu^{nr}}$-linear
nilpotent operator on $D(V)$
satisfying
$\sigma N=N{\frak p}^{n(\sigma)}N\sigma$
since
$\varphi N=p N\varphi$.
Thus an $\widehat{L_\mu^{nr}}$-linear
action ${}'\rho_{\mu,\pi,v}$
of the Weil-Deligne group
on $D(V)$ is defined.

We apply the above construction
$V\mapsto D(V)$
to the restriction
$\rho_{f,\mu,{\frak p}}$
of the $p$-adic representation
associated to $\pi_f$
to the decomposition group 
$\text{Gal}(\bar F_{\frak p}/F_{\frak p})$
for a place ${\frak p}|p$.
Thus we obtain an $\widehat{L(f)_\mu^{nr}}$-representation
${}'\rho_{f,\mu,{\frak p}}$
of the Weil-Deligne group
$'{}W(\bar F_{\frak p}/F_{\frak p})$.
Our main result is the following.

\proclaim{Theorem 1}
Let the assumptions be the same as in Theorem 0 above
and let $\mu$ be a place of $L(f)$ dividing
the characteristic of a prime ${\frak p}$ of $F$.
Then, the representation
$\rho_{f,\mu,{\frak p}}$ 
of $\text{Gal}(\bar F_{\frak p}/F_{\frak p})$
is potentially semi-stable
and there is an isomorphism
$${}'\rho_{f,\mu,{\frak p}}^{\text{F-ss}}\simeq
\check\sigma_h(\pi_{f,{\frak p}})$$
of representations of
the Weil-Deligne group
${}'W(\bar F_{\frak p}/F_{\frak p})$.
\endproclaim

\demo{Remark}
By the semi-stability of $\rho_{f,\mu,{\frak p}}$,
the representation ${}'\rho_{f,\mu,{\frak p}}$
is of degree 2.
Similarly as in the $\ell$-adic case,
Theorem implies that the left hand side
${}'\rho_{f,\mu,{\frak p}}^{\text{F-ss}}$
is $L(f)$-rational.
\enddemo

By the argument using a quadratic base change
as in [C2],
we may assume there exists a finite place
${\frak v}\neq {\frak p}$ where
$\pi_{f,{\frak v}}$ lies in the discrete series
in the case where $g=[F:\Bbb Q]$ is even.

We will prove Theorem 1 by comparing
$p$-adic cohomology with $\ell$-adic cohomology.
Let $\lambda$ be a place of $L(f)$ dividing 
a prime $\ell\neq p$.
By Theorem 0 applied to $\rho_{f,\lambda,{\frak p}}$,
it is enough to compare ${}'\rho_{f, \lambda,{\frak p}}$
with ${}'\rho_{f,\mu,{\frak p}}$.
More precisely, we prove the following.

\proclaim{Claim 1}
Let the notation be as in Theorem.
Let ${\frak p}|p$ be a finite place of $F$ and
let $\lambda$ and $\mu$ be places of $L(f)$
dividing $\ell\neq p$ and $p$ respectively.
Then the following holds.

\noindent
(0) The representation
$\rho_{f,\mu,{\frak p}}$ is potentially semi-stable.

\noindent
(1)
For $\sigma\in W^+=
\{\sigma\in W(\bar F_{\frak p}/F_{\frak p})|
n(\sigma)\ge0\}$, we have 
an equality in some finite extension of $L(f)$
$$\text{Tr }{}'\rho_{f,\lambda,{\frak p}}(\sigma)
=
\text{Tr }{}'\rho_{f,\mu,{\frak p}}(\sigma).
$$

\noindent
(2)
Let $N_\lambda$ and
$N_\mu$ be
the nilpotent monodromy operators
for ${}'\rho_{\lambda,\pi,{\frak p}}$ and
${}'\rho_{\mu,\pi,{\frak p}}$ respectively.
Then
$N_\lambda=0$ if and only if
$N_\mu=0.$
\endproclaim

By Lemma 1 [Sa],
Theorem 1 follows from Claim 1.
In (1), we may allow a finite extension
since we already know that the
left hand side is in $L(f)$.

The assertion (0) is a special case of (1)
where $\sigma=1$.
We deduce the assertion (2)
from (1) together with the
monodromy-weight conjecture,
Theorem 2 below,
asserting that the monodromy filtration
gives the weight filtration.

Let $V$ be a representation
of the Weil-Deligne group ${}'W_{\frak p}$.
We assume $N^2=0$.
Then $0\subset W_{-1}V=\text{Image }N\subset
W_0V= \text{Ker }N
\subset W_1V=V$
is a filtration by subrepresentations of $V$.
It is called the monodromy filtration.
We put $\text{Gr}^W_1(V)=V/\text{Ker }N,
\text{Gr}^W_0(V)=\text{Ker }N/\text{Image }N$ and
$\text{Gr}^W_{-1}(V)=\text{Image }N$.
Then each graded piece is a representation of
the Weil group.
The monodromy operator $N$
induces an isomorphism
$\text{Gr}^W_1(V)(1)\to
\text{Gr}^W_{-1}(V)$.
For a lifting $F$ of the geometric Frobenius
$Fr$, the eigenvalues upto roots of unity 
are independent of a choice of lifting.
We say an algebraic number is
pure of weight $n$ if
the complex absolute value
of its conjugates are $N\frak p^{\frac n2}$.
Then,
for an integer $n\in \Bbb Z$,
we say that the monodromy
filtration of $V$ is pure of
weight $n$,
if the eigenvalues of a lifting $F$
of $Fr$ acting on $Gr^W_i$ for each $i$
are algebraic numbers of weight $n+i$.

\proclaim{Theorem 2}
Let the notation be as in Claim 1.
Then the monodromy filtration
of the representations
${}'\rho_{f,\lambda,{\frak p}}(F)$
and ${}'\rho_{f,\mu,{\frak p}}(F)$
of the
Weil-Deligne group are pure of weight $w-1$.
In other words,
the eigenvalue $\alpha$
of ${}'\rho_{f,\lambda,{\frak p}}(F)$ for
an arbitrary lifting 
$F\in W(\bar F_{\frak p}/F_{\frak p})$
of the geometric Frobenius
is of weight $n$ where
$$n=
\cases 
w-1\quad&\text{if $N=0$}\\
w-2\quad&\text{if $N\neq0$ and $\alpha$
is the eigenvalue on Ker $N$}\\
w\quad&\text{if $N\neq0$ and $\alpha$
is the eigenvalue on Coker $N$.}
\endcases
$$
\endproclaim

\demo{Remark}
The assertion for the case $N\neq0$ is easy since
we know the determinant and 
$N:\text{Gr}^W_1(V)(1)\to
\text{Gr}^W_{-1}(V)$
is an isomorphism.
\enddemo

We show that Theorem 2 and assertion (1) in Claim 1
imply (2) in Claim 1.
In fact, by (1), the eigenvalues of a lifting
$F$ of Frobenius are the same for $\lambda$ and $\mu$.
By Theorem 2, we distinguish
the 2 cases $N=0$ and $N\neq 0$ by their absolute values.
Thus (2) follows from (1) and Theorem 2.

Thus Theorem 1 is reduced to
the assertion (1) in Claim 1
and Theorem 2.

\medskip
\subheading{3. Cohomological construction
of the $\ell$-adic representation}

Carayol constructs 
an $\ell$-adic representation
associated to a Hilbert modular form
by decomposing the etale cohomology
$H^1(M_{K,\bar F},\Cal F_\lambda)$
of a Shimura curve with a coefficient sheaf
$\Cal F_\lambda$.
Here, we briefly recall 
the construction with a slight modification.
Using the construction,
we give a statement, Claim 2,
which implies the main results.

First we recall the definition of the Shimura curve.
We choose and fix a real place $\tau_1$ of
the totally real field $F$
and regard $F$ as a subfield of $\Bbb R\subset \Bbb C$
by $\tau_1$.
When the degree $g=[F:\Bbb Q]$ is
even, we also fix a finite place $v_0$.
Let $B$ be a quarternion algebra over $F$
ramifying exactly at
the other real places $\{\tau_2,\ldots,\tau_g\}$
if $g=[F:\Bbb Q]$ is odd
and at $\{\tau_2,\ldots,\tau_g,\frak v\}$
if $g$ is even.

Let $G$
denote $B^\times$
regarded as an algebraic group
over $\Bbb Q$.
Let $X$ be the $G(\Bbb R)$-conjugacy of
the map
$$\matrix
h:&\Bbb C^\times&\to
&G(\Bbb R)=(B\otimes_{\Bbb Q}\Bbb R)^\times\simeq
&GL_2(\Bbb R)\times \Bbb H^\times
\cdots\times \Bbb H^\times,\\
&a+b\sqrt{-1}&\mapsto &&
\left(\pmatrix a&b\\-b&a\endpmatrix,1,\ldots,1\right).
\endmatrix$$
The conjugacy class $X$ is naturally identified with
the union $\Bbb P^1(\Bbb C)-\Bbb P^1(\Bbb R)$
of the upper and lower half planes.
Let $M=M(G,X)=(M_K)_K$
be the canonical model of Shimura variety
defined for $G$ and $X$.
It is defined over the reflex field $F$.
Here $K$ runs the
open compact subgroups of
$G(\Bbb A^f)=(B\otimes_F\Bbb A_{F,f})^\times$.
Each $M_K$ is a
proper smooth 
but not necessarily geometrically connected
curve over $F$.
Since the reciprocity map
$F^\times\to G^{ab}=F^\times$ is the identity,
the constant field $F_K$ of $M_K$
is the abelian extension of $F$ corresponding
to the compact open subgroup
$\text{Nrd}_{B/F}K\subset \Bbb A_{F,f}^\times$.
The projective system $(M_K)_K$
has a natural right action of the finite adeles
$G(\Bbb A_f)$.
For $g\in G(\Bbb A^f)$
and open compact subgroup $K,K'\subset G(\Bbb A^f)$
such that $g^{-1}Kg\subset K'$,
we have 
$g:M_K\to M_{K'}$.
The set of $\Bbb C$-valued points
$M_K(\Bbb C)$
are identified with the set of double cosets
$G(\Bbb Q)\backslash X\times G(\Bbb A^f)/K.$
The action of $G(\Bbb Q)=B^\times$ on $X$
is induced by $B^\times \to (B\otimes_{F,\tau_1}\Bbb R)^\times
\simeq GL_2(\Bbb R)$.
For $g,K,K'$ as above,
the map
$g:M_K(\Bbb C)\to M_{K'}(\Bbb C)$
is induced by $(x,g_1)\mapsto (x,g_1g).$

We will define a smooth $L_\lambda$-sheaf $\Cal F^{(k)}_\lambda$
on the Shimura curve $M$.
It is the dual of the sheaf 
denoted $\Cal F_\lambda$ in [Ca].
We prefer the dual because it is 
related directly to a direct summand of
a cohomology sheaf
as we will see in later sections.
Let $k=((k_1,\ldots,k_g),w)$ be a multiweight
and put $n=n(k)=\prod_i(k_i-1)$.
The algebraic group denoted
$G^c$ in [Mi] Chap. III
for our group $G=B^\times$ 
is the quotient of $G$ by $\text{Ker}(N_{F/\Bbb Q}:
F^\times \to \Bbb Q^\times)$.
Here we identify algebraic groups over $\Bbb Q$ and
their $\Bbb Q$-valued points and $F^\times\subset B^\times$
denotes the center of $G$.
In order to define the sheaf $\Cal F^{(k)}_\lambda$,
We define a representation of algebraic group
$\rho=\rho^{(k)}:G\to GL_n$
factoring the quotient $G^c$ as follows.
We have
$B\otimes_{\Bbb Q}\Bbb C \simeq M_2(\Bbb C)^I$
where $I=\{\tau_1,\ldots,\tau_g\}$
is the set of embeddings $F\to \Bbb C$.
It induces an isomorphism
$G_{\Bbb C}\to GL_{2,\Bbb C}^I$.
We define the morphism
$\rho=\rho^{(k)}:G\to GL_n$
to be the composite of this isomorphism
with the tensor product
$\bigotimes_{i\in I}((
\text{Sym}^{k_i-2}\otimes\det{}^{(w-k_i)/2})\circ \check{pr}_i)$.
Here $\check{pr}_i$ denotes the contragradient representation
of the $i$-th projection $GL_{2,\Bbb C}^I
\to GL_{2,\Bbb C}$.
Since the restriction of the center $F^\times$ is the
multiplication by $N_{F/\Bbb Q}^{-(w-2)}$,
it factors through the quotient $\rho^{(k)}:G^c\to GL_n$.

We take a number field $L\subset \Bbb C$ where
$\rho=\rho^{(k)}:G\to GL_n$ is defined.
For example
it is enough to take
$L$ containing the conjugates of $F$ and splitting $B$.
We identify $\{\tau_i:F\to L\}=
\{\tau_i:F\to \Bbb C\}$
by the inclusion $L\to \Bbb C$.
We define the smooth $L_\lambda$-sheaf
$\Cal F_\lambda^{(k)}$ on $M$
to be the $L_\lambda$-component of
the smooth $L\otimes\Bbb Q_\ell$-sheaf
$V_\ell(\rho^{(k)})$ attached to the representation $\rho^{(k)}$
(loc.\ cit.\ 6).
We consider the inductive limit
$$H^1(M_{\bar F},\Cal F_\lambda^{(k)})
=\varinjlim_K
H^1(M_{K,\bar F},\Cal F_{K,\lambda}^{(k)}).
$$
By the natural action of $G(\Bbb A_f)$
on the projective system
$(M_K,\Cal F_{K,\lambda}^{(k)})_K$,
it is a representation of
$G(\Bbb A_f)\times \text{Gal}(\bar F/F)$.
The structure as a birepresentation
is described as follows.

\proclaim{Lemma 1}
Let $k$ be a multiweight and $L\subset \Bbb C$ be
a number field splitting $F$ and $B$
as above.
Then, we have the following.

\noindent
(1).
Let $\pi_f$ be a cuspidal automorphic representation $\pi_f$
of $GL_2(\Bbb A_{F,f})$ of multiweight $k$.
If $g=[F:\Bbb Q]$ is even, we assume that
the $\frak v$-component $\pi_{f,\frak v}$
is in the discrete series.
Then the representation $\pi'_f$
of $G(\Bbb A_{f})$ corresponding to 
$\pi_f$ by the Jacquet-Langlands correspondence
has an $L(f)$-structure $\pi'_{f,L(f)}$.

\noindent
(2).
There exists an isomorphism 
$$H^1(M_{\bar F},\Cal F_\lambda^{(k)})
\simeq
\bigoplus_{f}
\left(\pi'_{f,L(f)}\otimes_{L(f)}
\bigoplus_{\lambda'|\lambda} V_{\lambda', f}
\right).$$
of representations of $G(\Bbb A_f)\times 
{\text{\rm Gal}}(\bar F/F)$.
Here $f$ runs the conjugacy classes over $L$
of normalized eigen newforms
of multiweight $k$,
such that,
when $g=[F:\Bbb Q]$ is even,
the $\frak v$-component $\pi_{f,\frak v}$
lies in the discrete series.
\endproclaim

Using Lemma 1, we reduce Theorems to a
statement below, Claim 2, on the cohomology.
Let $K\subset G(\Bbb A_f)$ be
a sufficiently small open compact subgroup.
We take a
sufficiently divisible integral ideal 
$\frak n$ of $O_F$,
divisible by $\frak p$ and by
$\frak v$ if $g$ is even.
We assume $K$ is of the form
$K=K_{\frak n}K^{\frak n}$.
Here
$K_{\frak n}\subset
\prod_{{\frak r}|{\frak n}}B_{\frak r}^\times$
is an open compact subgroup
and $K^{\frak n}=
\prod_{{\frak r}\nmid {\frak n}}GL_2(O_{F_{\frak r}})$
for some isomorphism
$\prod'_{{\frak r}\nmid {\frak n}}
B_{\frak r}
\simeq \prod'_{{\frak r}\nmid {\frak n}}M_2(F_{\frak r})$.
Let $T^{\frak n}=
L[T_{\frak r};\ {\frak r}\nmid {\frak n}]$ be 
the free $L$-algebra
generated by the Hecke operators
$T_{\frak r}$ for ${\frak r}\nmid {\frak n}$.
We consider
$H^1(M_{K,\bar F},\Cal F_\lambda^{(k)})$
as a $T^{\frak n}$-module.

\proclaim{Claim 2}
Let $K\subset G(\Bbb A_f)$
be a sufficiently small
open compact subgroup
and let $\frak n\subset O_F$ be
a sufficiently divisible ideal.
We assume $K=K_{\frak n}K^{\frak n}$ as above.
Then,

\noindent (0)
The representation 
$H^q(M_{K,\bar F},\Cal F_\lambda^{(k)})$
of $G_{F_\frak p}$
for $q=0,1,2$ 
is potentially semi-stable.

\noindent (1)
For $\sigma\in W^+$ and $T\in T^{\frak n}$,
we have equalities in a finite extension of $L$
$$
\sum_{q=0}^2(-1)^q\text{Tr}(\sigma\circ T|
 H^q(M_{K,\bar F},\Cal F_\lambda^{(k)}))
=
\sum_{q=0}^2(-1)^q
\text{Tr}(\sigma\circ T|
D(H^q(M_{K,\bar F},\Cal F_\mu^{(k)})))
$$

\noindent (2)
For the representations 
$H^1(M_{K,\bar F},\Cal F_\lambda^{(k)})$
and 
$D(H^1(M_{K,\bar F},\Cal F_\mu^{(k)})))$
of the Weil-Deligne
group ${}'W_{F_{\frak p}}$,
their monodromy filtrations 
are pure of weight $w-1$.
\endproclaim

We prove that the assertions in Claim 2
imply the corresponding 
assertions (0) and (1) in
Claim 1 and Theorem 2,
admitting Lemma 1.
Let $f$ be a normalized eigen new cuspform
of multiweight $k$.
By Lemma 1,
we have a cuspidal automorphic representation
$\pi'_f$ of $G'(\Bbb A_f)$ defined over $L(f)$.
Replacing $L$ by $L(f)$ if necessary,
we may assume $L=L(f)$.
Let $K$ be a sufficiently small open compact
subgroup satisfying
$\pi_f^{\prime K}\neq 0$ and Claim 2.
The representations
$V_{f,\lambda}$ and $V_{f,\mu}$
are direct summands of
$H^1(M_{K,\bar F},\Cal F_\lambda^{(k)})$
and $H^1(M_{K,\bar F},\Cal F_\mu^{(k)})$
by Lemma 1 respectively.
Hence the assertion (0) in Claim 2 implies
the assertion (0) in Claim 1
and the assertion (2) in Claim 2 implies
Theorem 2.

We show that the equality (1) of the traces
in Claim 2
implies the equality (1) in Claim 1.
First we show that
the equality for the alternating sum
implies the equality for each piece
$$\text{Tr}(\sigma\circ T|
H^r(M_{K,\bar E},\Cal F_\lambda^{(k)}))
=\text{Tr}(\sigma\circ T|
D(H^r(M_{K,\bar E},\Cal F_\mu^{(k)})))$$
for $r=0,1,2$.
In fact, it is sufficient to show
the equality for $r=0,2$.

We show that $H^0=H^2=0$
if $k\neq((2,\cdots,2),w)$.
The fundamental group $\pi_1(M_{K,\bar F})$
of the geometric fiber 
is isomorphic to $\text{Ker}
(\text{Nrd}_{B/F}:K\to \hat O_F^\times)$.
Hence its Lie algebra generates
$B^0=\text{Ker}
(\text{Trd}_{B/F}:B\to F)$ over $F$.
Hence the Lie algebra is
$B^0\otimes_{\Bbb Q}
L\simeq \frak s\frak l_2(L_\lambda)^g$
is also generated by the Lie algebra of
$\pi_1(M_{K,\bar F})$ over $L$.
It follows easily from this that
the representation of $\pi_1(M_{K,\bar F})$
corresponding to the sheaf $\Cal F_{\lambda}$
hence the sheaf itself is irreducible.
Hence its largest geometrically constant
subsheaf and quotient sheaves
are 0 unless $k=((2,\cdots,2),w)$.

We assume $k=((2,\cdots,2),w)$
and we show the equality for $r=0,2$.
Then the sheaf $\Cal F^{(k)}_\lambda$
is defined by the character 
$N_{F/\Bbb Q}^{-\frac{w-2}2}
\circ \text{Nrd}_{B/F}:G\to \Bbb G_m$
and is isomorphic to the Tate twist
$L_\lambda(-\frac{w-2}2)$.
It is sufficient to show
the assertion for $H^0$
since $H^2\simeq H^0(-1)$.
Let $F_K=\Gamma(M_K,O)$ be the constant field
of $M_K$.
Then
there is an isomorphism
$$H^0(M_{\bar F},\Cal F^{(k)}_\lambda)
\simeq
\varprojlim_K H^0(F_{K,\bar F},
L_\lambda(-\frac{w-2}2))$$
of $G_F\times G(\Bbb A_f)$-module.
On the right hand side,
the Galois action is the natural one.
The action of $G(\Bbb A_f)$
is defined by that induced by
its action on
$\varprojlim_K\text{Spec }F_K$
multiplied by the character
$$G(\Bbb A_f)@>{N_{F/\Bbb Q}^{\frac{w-2}2}
\circ \text{Nrd}_{B/F}}>>\Bbb A_f^\times/\Bbb Q^{+\times}
@<{\sim}<<\hat \Bbb Z^\times
@>>>\Bbb Z_\ell^\times
\subset L_\lambda^\times.$$
From this, we easily deduce the equality
for $r=0$.

We deduce the equality (1) in Claim 1
from the equality above for $r=1$.
By the strong multiplicity one theorem,
the image of the Hecke algebra $T^{\frak n}$
in $\text{End}_L(S_L^K)$
is $\prod_{f'}L(f')$ where $f'$
runs the conjugacy class of
eigen newforms $f'$
as in Lemma 1 such that ${\pi'_{f'}}^K\neq 0$.
Let $e\in T^{\frak n}$ be an element
whose image is the idempotent corresponding to
the component $L(f)=L$.
Then if we put $d=\dim \pi^{\prime K}_f$,
we see that $e\cdot H^1(M_{K,\bar F},
\Cal F_\lambda^{(k)})$
is isomorphic to the direct sum $\check\sigma_h(\pi)^{\oplus d}$
by Lemma 1 (2).
Hence we have
$$\align
d\cdot \text{Tr}'\rho_{\lambda,f,{\frak p}}(\sigma)
=&
\text{Tr}(\sigma\circ e|
 H^1(M_{K,\bar F},\Cal F_\lambda^{(k)}))
\\
d\cdot \text{Tr}'\rho_{\mu,f,{\frak p}}(\sigma)
=&
\text{Tr}(\sigma\circ e|
D( H^1(M_{K,\bar F},\Cal F_\mu)).
\endalign$$
Thus the equality (1) in Claim 1
follows from that in Claim 2.
It is clear that the assertion (2) in Claim 2
implies Theorem 2.

Therefore Theorems 1 and 2
are reduced to Claim 2 and Lemma 1. 

\demo{Proof of Lemma 1}
We will define an admissible representation $S'_L$
of $G(\Bbb A_f)$ over $L$ satisfying the following properties

(1')
$S'_L\otimes_L\Bbb C\simeq \bigoplus_f\pi'_f$ as
a representation of $G(\Bbb A_f)$ over $\Bbb C$.

 (2')
$H^1(M_{\bar F},\Cal F_\lambda^{(k)})
\simeq S_L^{\prime \oplus 2}\otimes_LL_\lambda$
as a representation of $G(\Bbb A_f)$ over $L_\lambda$.

\noindent
First we prove Lemma 1 assuming we have such a representation
$S'_L$.
We show $\pi'_f$ is defined over $L(f)$.
If $g=[F:\Bbb Q]$ is odd,
we have $G(\Bbb A_f)\simeq GL_2(\Bbb A_{F,f})$ and $\pi_f=\pi'_f$
and there is nothing to prove.
We show the case $g$ is even.
It is enough to show that each factor
$\pi'_{f,{\frak r}}$ of $\pi'_f=\bigotimes_v\pi'_{f,{\frak r}}$ is
defined over $L(f)$.
Let $\frak n$ be the level of $f$ and
$K_1({\frak n})=K_1({\frak n})_{\frak r}
\cdot K_1({\frak n})^{\frak r}\subset GL_2(\Bbb A_{F,f})$.
Then the representation $\pi_{f,{\frak r}}$
is given as the fixed subspace 
$\pi_{f,{\frak r}}=
\pi_f^{K_1({\frak n})^{\frak r}}$ and is defined over $L(f)$.
For ${\frak r}\neq {\frak q}$,
we have $\pi_{f,{\frak r}}=\pi'_{f,{\frak r}}$ and it is
defined over $L(f)$.
Finally we consider the case ${\frak r}=\frak q$.
Then by (1), we see that
the intertwining space
$\text{Hom}_{G(\Bbb A_f^{\frak q})}
(\bigotimes_{{\frak r}\neq \frak q}\pi'_{f,\frak r,L(f)},
S_L\otimes_LL(f))$ is an $L(f)$-structure of $\pi'_{f,\frak q}$.

Next we show the isomorphism (2).
We put
$$V_{\lambda',f}=\text{Hom}_{G(\Bbb A_f)}
\left(\pi'_{f,L(f)}\otimes_{L(f)}L(f)_{\lambda'},
H^1(M_{\bar F},\Cal F_\lambda^{(k)})\right).$$
Then by (2'), each $V_{\lambda',f}$
is an $L(f)_{\lambda'}$-representation of
$G_F$ of degree 2 and we have
$$H^1(M_{\bar F},\Cal F_\lambda^{(k)})
\simeq
\bigoplus_{f}
\left(\pi'_f\otimes_{L(f)}
\bigoplus_{\lambda'|\lambda} V_{\lambda', f}
\right).$$
Therefore it is sufficint to show that
the $\ell$-adic representation
$V_{\lambda', f}$ is associated to $f$:
$V_{\lambda', f}\simeq\check \sigma_h(\pi_f)$.
We may extend the scaler to $\bar L_\lambda$.
Hence it is enough to show that
$H^1(M_{\bar F},\Cal F_\lambda^{(k)})\otimes \bar L_\lambda$
admits a direct sum decomposition of the form
$$H^1(M_{\bar F},\Cal F^{(k)}_\lambda)\otimes \bar L_\lambda
\simeq
\bigoplus_\pi\pi\otimes \check\sigma_h(\pi).$$
In [C2], it is shown that for $\Cal F_\lambda$,
we have a direct sum decomposition of the form
$H^1(M_{\bar F},\Cal F_\lambda)\simeq
\bigoplus_{\pi'}\pi'\otimes \check\sigma_h(\pi')$.
Since $\Cal F^{(k)}_\lambda$ here is the dual of
$\Cal F_\lambda$ there,
we have
$H^1(M_{\bar F},\Cal F^{(k)}_\lambda)\simeq
\bigoplus_{\pi'}\check\pi'\otimes \sigma_h(\pi')(-1)$
by Poincar\'e duality.
Since $\check \sigma_h(\check \pi')\simeq
\sigma_h(\pi')(-1)$,
the claim is proved.

We define the space $S'_L$.
First we will define the automorphic vector bundle ([Mi] Chap.III)
$\Cal V(\Cal J)$ associated to
a $G^c$-equivariant vector bundle
$\Cal J=\Cal J^{(k)}$
on the compact dual $\check X$ 
and its canonical model $\Cal V(\Cal J)_L$.
Then it will be defined as
the limit of the spaces of global sections
$$S'_L=
\Gamma(M\otimes_FL,\Omega^1_M\otimes \Cal V(\Cal J)_L)=
\varinjlim_K
\Gamma(M_K\otimes_FL,\Omega^1_M\otimes \Cal V(\Cal J)_L).$$
We use the notation loc.\ cit.
The compact dual $\check X$ is $\Bbb P^1_{\Bbb C}$
in our case. 
We define a $G^c$-equivariant vector bundle $\Cal J=\Cal J^{(k)}$
on $\check X$ in the following way.
Let $\omega$ be the dual of the tautological quotient bundle
on $\check X=\Bbb P^1_{\Bbb C}$.
We put
$\Cal J^{(k)}=
\omega^{\otimes k_1-2}\otimes\bigotimes_{i=2}^g
\text{Sym}^{k_i-2}(\Bbb C^{\oplus 2}).$
We define the action of $G_{\Bbb C}={GL_{2,\Bbb C}}^I$
on $\Cal J^{(k)}$ by giving the action of each factor
in the following way.
The first factor $GL_{2,\Bbb C}$ acts on $\check X$
in the natural way.
On $\omega^{\otimes k_1-2}$,
we consider $\det^{-\frac{w-k_1}2}$-times the natural action.
For $i\neq 1$,
the $i$-th factor $GL_{2,\Bbb C}$ acts on $\check X$ trivially.
On $\text{Sym}^{k_i-2}(\Bbb C^{\oplus 2})$,
we consider 
$\det^{-\frac{w-k_i}2}$-times the action induced by the 
contragradient action of $GL_2$.
By taking the tensor product,
we obtain a $G_{\Bbb C}$-equivariant bundle $\Cal J=\Cal J^{(k)}$.
Since the center $\Bbb G_m^I$ acts by $-(w-2)$-nd power
of the product character,
it defines a $G^c_{\Bbb C}$-equivariant bundle.
It is clearly defined over the number field $L\supset F$.
Hence by [Mi] Chap.III Theorem 5.1.(a),
we obtain a $G(\Bbb A_f)$-equivariant vector bundle
$\Cal V(\Cal J)_L$ on $M_L$.
Thus the representation
$S'_L=
\Gamma(M\otimes L,\Omega^1_M\otimes \Cal V(\Cal J)_L)$
is defined.

By the Jacquet-Langlands correspondence [JL],
we have $S'_{\Bbb C}=S'_L\otimes_L\Bbb C\simeq\bigoplus_f\pi'_f$
where $f$ runs cuspidal automorphic representation
of $GL_2(\Bbb A_{F,f})$ of multiweight $k$
such that, if $[F:\Bbb Q]$ is even,
the $q$-component $\pi_{f,q}$ is in the discrete series.
Hence the representation
$S'_L$ satisfies the property (1') above.

We show the isomorphism (2').
Attached to the representation
$\rho^{(k)}:G^c\to GL_n$ defined over $L$,
we have a local system
$\Cal F^{(k)}= V(\rho^{(k)})$
of $L$-vector spaces on $M(\Bbb C)$.
Since 
$H^1(M_{\bar F},\Cal F_\lambda^{(k)})\simeq
H^1(M(\bar \Bbb C),\Cal F^{(k)})\otimes_LL_\lambda$,
it is sufficient to show
$H^1(M(\bar \Bbb C),\Cal F^{(k)})\otimes_L\Bbb C
\simeq S_L^{\prime \ \oplus 2}\otimes \Bbb C$
as a representation of $G(\Bbb A_f)$.
We deduce it from Hodge decomposition,
which is a generalization of the Eichler-Shimura isomorphism.
Let $\Cal F^{(k)}_{\Bbb C}=\Cal F^{(k)}\otimes_L\Bbb C$.
We regard it as a local system of $\Bbb R$-vector spaces
endowed with a ring homomorphism 
$\Bbb C\to \text{End}(\Cal F^{(k)}_{\Bbb C})$.
We consider the filtration on $\Cal F^{(k)}_{\Bbb C}
\otimes_{\Bbb R}O_{M(\Bbb C)}$ defined by $\rho\circ h_x$.
It defines on $\Cal F^{(k)}_{\Bbb C}$
a structure of variation of polarizable $\Bbb R$-Hodge
structures of weight $w-2$.

We put $\Cal F^{(k)}_{\Bbb C}
\otimes_{\Bbb C}O_{M(\Bbb C)}=\Cal V^{(k)}
(=\Cal V(\Cal \rho^{(k)}))$ and
let $\sigma:M(\Bbb C)\to M(\Bbb C)$ denote the complex conjugate.
We identify
$\Cal F^{(k)}_{\Bbb C}
\otimes_{\Bbb R}O_{M(\Bbb C)}= \Cal V^{(k)}\oplus \sigma^*\Cal V^{(k)}$.
The Hodge filtration $F^{w-2}(\Cal V^{(k)}\oplus \sigma^*\Cal V^{(k)})$
is given by $\Cal V(\Cal J^{(k)})\oplus \sigma^*\Cal V(\Cal J^{(k)})$.
Hence the Hodge decomposition gives
a $G(\Bbb A_f)$-equivariant isomorphism
$$\align &\ H^1(M(\Bbb C),\Cal F^{(k)}_{\Bbb C})\simeq
H^1(M(\Bbb C),\Omega^\bullet_M\otimes \Cal V^{(k)})\\
\simeq &\
H^0(M(\Bbb C),\Omega^1_M\otimes \Cal V(\Cal J^{(k)}))\oplus
\sigma^*H^0(M(\Bbb C),\Omega^1_M\otimes \Cal V(\Cal J^{(k)})).
\endalign$$
Since
$S'_{\Bbb C}=H^0(M(\Bbb C),\Omega^1_M\otimes \Cal V(\Cal J^{(k)}))$
and its complex conjugate
$\sigma^*H^0(M(\Bbb C),\Omega^1_M\otimes \Cal V(\Cal J^{(k)}))$
is identified with
$H^0(M(\Bbb C),\Omega^1_M\otimes \Cal V(\sigma^*\Cal J^{(k)}))$,
it is enough to show that the 
$G(\Bbb A_f)$-equivariant bundle
$\Cal J^{(k)}$ on $\check X$
is isomorphic to 
its complex conjugate $\sigma^*\Cal J^{(k)}$.
It follows immediately from that
the $GL_2$-action on the tautological quotient bundle
on $\Bbb P^1$
is defined over $\Bbb R$ and
the standard representation $\Bbb H^\times \to GL_2$
defined over $\Bbb C$
is $GL_2(\Bbb C)$-conjugate to its complex conjugate.
\enddemo

\medskip
\subheading{4. Shimura curves and sheaves on them}

We give a geometric construction of 
a certain pull-back of the
sheaf $\Cal F^{(k)}$
using functoriality of Shimura varieties.

First, we recall the definition of several
Shimura varieties introduced in [C1].
We take an imaginary quadratic field
$E_0=\Bbb Q(\sqrt{-a})$.
We fix an embedding $E_0\subset \Bbb C$.
We assume that the prime $p$ splits in $E_0$.
We put $E=FE_0=F\otimes_{\Bbb Q} E_0$
and $D=B\otimes_F E=B\otimes_{\Bbb Q} E_0$.
We consider the reductive group
$G''=B^\times\times_{F^\times} E^\times
\simeq B^\times\cdot E^\times\subset D^\times$.
Here and in the following,
we use its $\Bbb Q$-valued points $G(\Bbb Q)$
to describe an algebraic group $G$ over $\Bbb Q$.
As in [C1],
the notation $B^\times\times_{F^\times} E^\times$
does not mean the fiber product but the amalgamate sum.
Let $G'$ be the inverse image of 
$\Bbb Q^\times\subset F^\times$
by the map $\nu=\text{Nrd}_{B/F}\times
N_{E/F}:G''\to F^\times$.
We also consider tori
$T=E^\times$ and $T_0=E_0^\times$.
We consider the
$G'(\Bbb R)$-conjugacy class
$X'$ 
(resp. $G''(\Bbb R)$-conjugacy class
$X''$) of the morphism
$$\matrix
h':&\Bbb C^\times&\to&
G'(\Bbb R)\subset
G''(\Bbb R)=&GL_2(\Bbb R)\cdot\Bbb C^\times\times
\Bbb H^\times\cdot \Bbb C^\times\times\cdots
\Bbb H^\times\cdot \Bbb C^\times,\\
&z&\mapsto& &\left(\pmatrix a& b\\ -b& a\endpmatrix
\otimes 1,1\otimes z,\ldots,1\otimes z\right).
\endmatrix$$
We also consider morphisms
$$\alignat 2
h_E:&\ \Bbb C^\times\to
T(\Bbb R)=
\Bbb C^\times\times\Bbb C^\times \times\cdots
\Bbb C^\times,\quad &
z&\ \mapsto (z^{-1},1,\ldots,1),\\
h_0:&\ \Bbb C^\times\to
T_0(\Bbb R)=
\Bbb C^\times &
z&\ \mapsto z^{-1}.
\endalignat$$
The conjugacy classes $X',X''$ have natural structures
of complex manifold and
are isomorphic to the upper half plane
$X^+$ and to the union of upper and lower half planes
$X$ respectively.
Let $M'=M(G',X'),M''=M(G'',X''),
N=M(T,h_E)$ and $N_0=M(T_0,h_0)$
be the canonical models of the
Shimura varieties
defined over the reflex fields
$E,E,E$ and $E_0$ respectively.
The reciprocity map
$E^\times\to E^\times$
is the identity for $(T,h_E)$.
For an open compact subgroup
$K\subset \Bbb A_{E,f}^\times$,
the canonical model $N_K$ is the spectrum
of the abelian extension $E_K$ 
corresponding to $K$
by class field theory.
The same thing holds for the canonical model of $N_0$.

We define morphisms between Shimura curves.
We consider the morphism
$\alpha:G\times T\to G''$ of algebraic groups
inducing
$$B^\times \times E^\times \to G''(\Bbb Q)
\subset (B\otimes E)^\times:\
(b, e)\to
b\otimes N_{E/E_0}(e)\cdot e^{-1}$$
on $\Bbb Q$-valued points.
Since $h'=\alpha \circ (h \times h_E)$,
it induces a homomorphism of Shimura varieties
$M\times N\to M''$
defined over $E$.
We let $\alpha$ also denote the
morphism $M\times N\to M''$.
The inclusion $G'\to G''$ induces
a natural map $M'\to M''$
of Shimura varieties over $E$.
Let $\beta:G\times T\to T_0$
be the morphism inducing
$N_{E/E_0}\circ \text{\rm pr}_2:
B^\times\times E^\times\to E_0^\times$ on
$\Bbb Q$-valued points.
Since $h_0=N_{E/E_0}\circ h$,
a homomorphism of Shimura varieties
$M\times N\to N_0$
defined over $E$ is induced.
We also let the map
$M\times N\to N_0$
be denoted by $\beta$.
We consider the diagram
$$\CD
M@<{pr_1}<< M\times N @>{\alpha}>> M'' @<<< M'\\
@.@V{\beta}VV\\
@.N_{0}
\endCD$$
of (weakly) canonical models of Shimura varieties over $E$.

We define an $L_\lambda$-sheaf
$\Cal F^{\prime\prime(k)}_\lambda$ on $M''$
analogous to $\Cal F^{(k)}_\lambda$.
Let $k=((k_1,\ldots,k_g),w)$ be the multiweight
and put $n=n(k)=\prod_i(k_i-1)$.
The algebraic group denoted
$G^{\prime\prime c}$ in [Mi] Chap. III
for the group $G''$
is the quotient of $G''$ by $\text{Ker}(N_{F/\Bbb Q}:
F^\times \to \Bbb Q^\times)$.
Here $F^\times$ is regarded as a
subgroup of the center $Z(G'')=E^\times$.
We define a representation of algebraic group
$\rho=\rho^{\prime\prime (k)}:G''\to GL_n$
factoring the quotient $G^{\prime \prime c}$ as follows.
Recall that we have an isomorphism
$B\otimes_{\Bbb Q}\Bbb C \simeq M_2(\Bbb C)^I$.
It induces an injection
$G''_{\Bbb C}\to (GL_{2,\Bbb C}\times GL_{2,\Bbb C})^I$.
For each $i\in I$, the first component
corresponds to the inclusion $E_0\to \Bbb C$
and the second one corresponds to its complex conjugate.
We define the morphism
$\rho''=\rho^{\prime\prime (k)}:G\to GL_n$
to be the composite of the injection
with the tensor product
$\bigotimes_{i\in I}((
\text{Sym}^{k_i-2}\otimes\det{}^{(w-k_i)/2})\circ \check{pr}_{2,i})$.
Here $\check{pr}_{2,i}$ denotes the contragradient representation
of the $(2,i)$-th projection $(GL_{2,\Bbb C}\times GL_{2,\Bbb C})^I
\to GL_{2,\Bbb C}$.
Since the restriction to 
the subgroup $F^\times\subset G''$
is the scalar
multiplication by $N_{F/\Bbb Q}^{-(w-2)}$,
it factors through the quotient 
$\rho^{\prime\prime (k)}:G^{\prime\prime c}\to GL_n$.
The morphism
$\rho''=\rho^{\prime\prime (k)}:G\to GL_n$ is defined over 
the composite field $LE_0$.
Replacing $L$ by $LE_0$ if necessary,
we assume it is defined over $L$.

We may also define it as follows.
Let $p_2:G''\to G$ be the map defined over $E_0$
induced by the second projection on
$(D\otimes_{\Bbb Q}E_0)^\times=
D^\times \times D^\times$
corresponding to the conjugate $E_0\to E_0$.
Then we have $\rho^{\prime\prime(k)}=
\rho^{(k)}\circ p_2$.

We define the smooth $L_\lambda$-sheaf
$\Cal F_\lambda^{\prime\prime (k)}$ on $M''$
to be the $L_\lambda$-component of
the smooth $L\otimes\Bbb Q_\ell$-sheaf
$V_\ell(\rho^{\prime\prime(k)})$ attached to the representation 
$\rho^{\prime\prime(k)}$
(loc.\ cit.\ 6).
By restriction, 
we obtain a smooth $L_\lambda$-sheaf
$\Cal F_\lambda^{\prime (k)}$ on $M'$
attached to the representation 
$\rho^{\prime(k)}=
\rho^{\prime\prime(k)}|_{G'}$.

We also define a sheaf $\Cal F(\chi)_\lambda$ on $N_0$.
The algebraic group $T_0^c$ in [Mi] Chap. III
is $T_0$ itself.
We define a character
$\chi:T_0\to \Bbb G_m$.
Over $\Bbb C$, we have
$T_{0,\Bbb C}\simeq \Bbb G_m\times \Bbb G_m$.
Here the first component
corresponds to the inclusion $E_0\to \Bbb C$
and the second one corresponds to its complex conjugate.
We define the morphism $\chi:T_0\to \Bbb G_m$
to be the inverse of the first projection.
We also define the morphism $\bar \chi$
to be the inverse of the second projection.
Their product $\chi_0=\chi\bar \chi$
is the inverse of the norm map 
$\chi_0=N_{E/\Bbb Q}^{-1}:T_0\to \Bbb G_m$.
They are defined over $E_0\subset L$.
We define the smooth $L_\lambda$-sheaf
$\Cal F(\chi)$ on $N_0$
to be the $L_\lambda$-component of
the smooth $L\otimes\Bbb Q_\ell$-sheaf
$V_\ell(\chi)$ attached to the representation $\chi$.
The sheaf
$\Cal F(\chi_0)$ is defined similarly.

We have $\rho^{''(k)}\circ \alpha=
(\rho^{(k)}\circ pr_1)\times 
(\bar \chi^{(w-2)(g-1)}\circ N_{E/E_0}\circ pr_2)$.
In other words, we have
a commutative diagram
$$\CD
G \times T @>{\alpha\times \beta}>> 
G^{\prime\prime}\times T_0\\
@V{\rho^{(k)}\circ pr_1}VV 
@VV{\rho^{\prime\prime (k)}\times
\chi^{(g-1)(w-2)} \chi_0^{-(g-1)(w-2)}}V \\
GL_n  @<<{\text{product}}< GL_n\times \Bbb G_m
\endCD$$
of homomorphisms defined over $L$.
By the commutativity of the diagram,
we obtain 
an isomorphism of smooth $L_{\lambda}$-sheaves
$$
pr_1^*\Cal F^{(k)}
\simeq
\alpha^*\Cal F^{\prime\prime(k)}
\otimes
\beta^*\Cal F(\chi)^{\otimes(g-1)(w-2)}
\otimes
\beta^*\Cal F(\chi_0)^{-\otimes(g-1)(w-2)}$$
on $M\times N$.
The isomorphism is equivariant
with respect to the action of
$G(\Bbb A_f)\times T(\Bbb A_f)$.

The sheaf
$\beta^*\Cal F(\chi_0)$
together with the action of $T(\Bbb A_f)$ on it
is identified as follows.
Let $\beta_1:N\to N_0$
denote the map induced by $N_{E/E_0}$.
It is sufficient to describe
$\beta_1^*\Cal F(\chi_0)$.
If we forget the action,
it is just the Tate twist $L_\lambda(-1)$.
The action of $T(\Bbb A_f)$ is that
induced by the natural action of
$T(\Bbb A_f)$ on $N$
multiplied by the character 
$$T(\Bbb A_f)@>{N_{E/\Bbb Q}}>>
\Bbb A_f^\times/\Bbb Q^{+\times}
@<{\sim}<<\hat \Bbb Z^\times
@>>>\Bbb Z_\ell^\times
\subset L_\lambda^\times.$$

Thus the geometric construction of
$pr_1^*\Cal F^{(k)}$ is reduced to that of
$\Cal F^{\prime\prime(k)}$ and that of
$\Cal F(\bar \chi)$.

Before constructing $\Cal F^{\prime\prime(k)}$ geometrically,
we will study
its restriction $\Cal F^{\prime(k)}$ to $M'$.
We prepare some notations.
We consider the representation
$$\rho':G'\subset G''\subset D^\times
@>{b\mapsto \bar b^{-1}}>>
D^\times \subset GL(D)$$
defined over $\Bbb Q$.
Since
the algebraic group $G^{\prime c}$ 
for $G'$ is equal to $G'$ itself,
the representation $\rho'$
gives rise to
a smooth $\ell$-adic sheaf
$\Cal F'_\ell$ on $M'$ for each prime $\ell$.
It is a smooth sheaf of 
$D\otimes_{\Bbb Q}\Bbb Q_\ell$-modules
of rank 1.

Recall that we have an isomorphism
$D\otimes_{\Bbb Q}L\simeq (M_2(L)\times M_2(L))^I$.
For each $i\in I$,
the first component corresponds
to the embedding $E_0\subset L\subset \Bbb C$
and the second to its conjugate.
For each $i\in I$,
let $e_i\in D\otimes_{\Bbb Q}L$ denote the idempotent
whose $(2,i)$-th component is $\pmatrix 1& 0\\ 0&0\endpmatrix$
and the other components are 0
under the isomorphism above.
For each finite place $\lambda|\ell$,
we regard the $L_\lambda$-sheaf
$\Cal F'\otimes_{\Bbb Q_\ell}L_\lambda$
as a
$D\otimes_{\Bbb Q}L\simeq (M_2(L)\times M_2(L))^I$-module.
For each $i\in I$,
let $\Cal F'_i$ denote the
$e_i$-part $e_i(\Cal F'\otimes_{\Bbb Q_\ell}L_\lambda)$.
it is easy to see that
$$\Cal F^{\prime (k)}=
\bigotimes_{i\in I}
(\text{Sym}^{k_i-2}\Cal F'_i\otimes
 ({\det}\Cal F'_i)^{\otimes\frac{w-k_i}2})$$
as a smooth $L_\lambda$-sheaf 
on $M'$ with an action of $G'(\Bbb A_f)$.

In section 6.1,
we will construct 
the sheaf $\Cal F'$ and the idempotents
$e_i$ after recalling a modular interpretation
of $M'$ in Section 5.
We will also construct
$\Cal F(\chi)$ on $N_0$ in a similar way.
After that, we study the relation between
$M'$ and $M''$
and extend $\Cal F'$ to $M''$
in section 6.2.

\medskip
\subheading{5. Modular interpretation of $M'$ and $N_0$}

We recall the modular interpretation of 
the Shimura curve $M'$ on the category of schemes over $E$ [C1] 2.3.
In the notation of [D2] (4.9) and (4.13),
we put $L=V=D$,
let the involution $*$ on $D=B\otimes_F E$ to
be the tensor product of the main involution of $B$ and
the conjugate of $E$ and
let $\psi$ be the non-generate alternating form on $D$
defined by
$$\psi(x,y)=\text{Tr}_{E/\Bbb Q}
(\sqrt{-a} \text{Trd}_{D/E}xy^*).$$
Then the group $G$ in (4.9) loc.\ cit.\ is $G'$ here and
$G_1$ in (4.13) is $G''$ here.

We prepare some terminology to formulate
a moduli problem for $M'$.
Let $O_D$ be a maximal order in $D$
stable under the involution $*$.
An abelian scheme $A$ over a scheme $S$ 
is called an $O_D$-abelian scheme over $S$
when a ring homomorphism
$m:O_D\to \text{End}(A)$ is given.
When $S$ is a scheme over $\text{Spec }E$,
for an $O_D$-abelian scheme $A$ on $S$,
we define direct summands
$\text{Lie}^2 A\supset \text{Lie}^{1,2} A$
of the $O_D\otimes_{\Bbb Z}O_S=
D\otimes_{\Bbb Q}O_S$-module
$\text{Lie} A$ as follows.
The submodule 
$\text{Lie}^2 A$
is defined to be the submodule on which
the action of $E_0\subset D$
and that of $E_0\subset O_S$
are the conjugate to each other 
over $\Bbb Q$.
Similarly
$\text{Lie}^{1,2} A$
is the submodule where
the action of $E\subset D$
and that of $E\subset O_S$
are the conjugate to each other
over $F$.
They are the same as the tensor products
$\text{Lie}^2 A=\text{Lie} A\otimes_{E_0\otimes E_0}E_0,
\text{Lie}^{1,2} A=\text{Lie} A\otimes_{E\otimes E}E$
and hence are direct summands.
If $A$ is an $O_D$-abelian scheme,
the dual $A^*$ is considered as an $O_D$-abelian scheme
by the composite map 
$m^*:O_D@>*>> O_D^{opp}@>m>> \text{End}(A)^{opp}
@>{*}>> \text{End}(A^*)$
where $opp$ denotes the opposite ring.
A polarization $\theta\in
\text{Hom}(A,A^*)^{sym}$ of
an $O_D$-abelian scheme $A$
is called an $O_D$-polarization if it is 
$O_D$-linear.

Let $K\subset \hat O_D^\times\subset
G'(\Bbb A_f)$ be a sufficiently small
compact open subgroup.
Take a maximal order $O_D$ of
$D$ and 
let $\hat O_D=O_D\otimes \hat \Bbb Z\subset 
D\otimes \Bbb A_f$ be
the corresponding maximal order.
We assume $K\subset  \hat O_D^\times.$
Let 
$\hat T\subset D\otimes \Bbb A_f$
be a $\hat O_D$-lattice
satisfying $\psi(\hat T,\hat T)
\subset \hat\Bbb Z$.
We define a functor $M'_{K',\hat T}$
on the category of schemes over $E$
as follows.
For a scheme $S$ over $E$ let
$M'_{K',\hat T}(S)$ be
the set of isomorphism class of 
the triples $(A,\theta,\bar k)$
where 

(1)
$A$ is an $O_D$-abelian scheme on $S$ of dimension 4g
such that $\text{Lie}^2 A=\text{Lie}^{1,2} A$
and that it is a locally free $O_S$-module of rank 2.

(2)
$\theta\in \text{Hom}(A,A^*)^{sym}$ 
is an $O_D$-polarization of $A$.

(3)
$\bar k$ is a $K$-equivalent class of
a $O_D\otimes \hat\Bbb Z$-linear
isomorphism $k:\hat T(A)\to \hat T$
such that there exists a $\hat\Bbb Z$-linear
isomorphism $k'$
making the diagram 
$$\matrix
\hat T(A)\times \hat T(A)&@>{(1,\theta_*)}>> 
\hat T(A)\times \hat T(A^*)@>>>& \hat\Bbb Z(1)\\
@V{k\times k}VV @VV{k'}V\\
\hat T\times \hat T& @>>{\text{Tr} \psi}>& \hat\Bbb Z
\endmatrix$$
commutative.

\noindent
It is shown in [C1] (2.3) (2.6.2) that the scheme
$M'_{K'}$ represents the functor $M'_{K',\hat T}$.
It is easily checked that the functor is 
independent of a choice of $\hat T$
upto uniquely defined canonical isomorphism.
Let $A_{K',\hat T}$ denote the universal abelian scheme
over $M_{K'}$.
They form a projective system
$A=(A_{K',\hat T})_{K',\hat T}$.

We give a modular interpretation of
the action of $G'(\Bbb A_f)$ on $M'$ and on $A$.
Let $g\in G'(\Bbb A_f)$ and $K,K'\subset G'(\Bbb A_f)$ be
sufficiently small open subgroups
satisfying $g^{-1}Kg\subset K'$.
We take a maximal order $O_D$
and let $\hat T$ and $\hat T'$ be
a $K$-stable $O_D\otimes \hat \Bbb Z$-lattice and
a $K'$-stable $O_D^g\otimes \hat \Bbb Z$-lattice of $V\otimes \Bbb A_f$
satisfying $g^{-1}\hat T\subset \hat T'$
and $\psi(\hat T,\hat T),\psi(\hat T',\hat T')\subset \hat \Bbb Z$.
The functor
$$g_*:\Cal M_K@>>> \Cal M_{K'},
\quad [(A,\theta,k)]\mapsto [(A',\theta',k')]$$
is described as follows.
(Ind-)Etale locally on $S$,
we take an isomorphism
$\tilde k:\hat T\to \hat T(A)$ in the $K$-equivalent class $k$
and identify $\hat T(A)$ with $\hat T$ by $\tilde k$.
Let $g_*:A\to A'$ be the isogeny of $O_D$-abelian schemes
such that $\hat T(A')=g \hat T'\supset \hat T=\hat T(A)$.
The $K'$-equivalent class $k'$
is the class of the isomorphism
$g:\hat T'\to g\hat T'=\hat T(A')$.
The pair $(A',k')$ is independent of the choice of
$\tilde k$.
The polarization $\theta'$ on $A'$
is the map making the diagram
$$\CD
A@>{\nu^+(g)\theta}>> A^*\\
@V{g_*}VV @AA{{}^tg_*}A\\
A'@>>{\theta'}> A^{\prime *}
\endCD$$
commutative.
Here $\nu^+:G'(\Bbb A_f)\to \Bbb Q^{\times +}$ is 
the composite
$G'(\Bbb A_f)@>{\nu}>> \Bbb A_f^\times\to 
\Bbb A_f^\times/\hat\Bbb Z^\times @<{\sim}<< \Bbb Q^{\times +}$.
We have the universal $O_D$-isogeny $g:A_{K,\hat T}\to g^*A_{K',\hat T'}$
and a commutative diagram
$$\CD
A_{K,\hat T}@>{g_*}>> A_{K',\hat T'}\\
@VVV @VVV\\
M'_{K}@>{g}>> M'_{K'}.
\endCD$$

For later use, we will extend the action of
$G'(\Bbb A_f)$ on $M'$ and on $A$ to a larger group
$\tilde G$.
Let $G''(\Bbb R)_+$ be the inverse image of
$GL_2(\Bbb R)^+\Bbb C^\times\subset 
GL_2(\Bbb R)\Bbb C^\times$
by the first projection
$G''(\Bbb R)\to 
GL_2(\Bbb R)\cdot \Bbb C^\times$
and let
$G''(\Bbb Q)_+=G''(\Bbb Q)\cap G''(\Bbb R)_+=
\{\gamma\in G''(\Bbb Q)|
\nu(\gamma)\text{ is totally positive}\}$.
We put $\tilde G=G''(\Bbb Q)_+\cdot G'(\Bbb A_f)
\subset  G''(\Bbb A_f)$.
We extend the action of $G'(\Bbb A_f)$ on $M'$
to an action of $\tilde G$.
For $g\in G''(\Bbb A_f)$
and open compact subgroups
$K'\subset G'(\Bbb A_f)$ and
$K''\subset G''(\Bbb A_f)$ such that
$g^{-1}K'g\subset K''$,
let $g:M'_{K'}\to M''_{K''}$
denote the composite
$M'_{K'}\to M''_{gK''g^{-1}}@>{g}>> M''_{K''}$.
For $g\in \tilde G$
and open compact subgroups
$K'_1,K'_2\subset G'(\Bbb A_f)$ such that
$g^{-1}K'_1g\subset K'_2$,
the map $g:M'_{K'_1}\to M'_{K'_2}$
is defined as follows.
We may take an open compact subgroup
$K''\supset K'_2$ of $G''(\Bbb A_f)$
such that the canonical map
$M'_{K'_2}\to M''_{K''}$
is an open immersion
(see Lemma 2 in section 6.2).
Then since
$M'_{K'}(\Bbb C)=
G'(\Bbb Q)\backslash  G'(\Bbb A_f)\times X'/K'=
G''(\Bbb Q)_+\backslash  \tilde G \times X'/K'$,
the image of $g:M'_{K'_1}\to M''_{K''}$ is contained in
$M'_{K'_2}$.
Hence the required map
$M'_{K'_1}\to M'_{K'_2}$ is induced.
The modular interpretation of the action of
$\tilde G$ on $M'$ is described in the same way as above.
The only modification is that
$\nu^+$ is extended to
$\tilde G$
as the composite
$\tilde G@>{\nu}>> F^{\times +}\Bbb A_f^\times\to 
F^{\times +}\Bbb A_f^\times/\hat\Bbb Z^\times @<{\sim}<< F^{\times +}$.

Similarly, we have a modular interpretation for
$N_0$ in terms of elliptic curves 
with complex multiplication by
$O_{E_0}$.
Let $H\subset \hat O_{E_0}^\times$
be a sufficiently small open subgroup. 
We take a fractional ideal $R\subset E_0$
satisfying $\text{Tr}_{E_0/\Bbb Q}(\sqrt a R\bar R)
\subset \Bbb Z$.
Let
$\hat R=R\otimes \hat O_E$ be
the corresponding ideal.
We define a functor
$N_{0,H,\hat R}$
on the category of schemes over $E_0$
as follows.
For a scheme $S$ over $E_0$,
let
$N_{0,H,\hat R}(S)$ be the set of
isomorphism of the pairs
$(A,\bar k)$
where 

(1)
$A$ is an elliptic curve
endowed with a ring homomorphism
$O_{E_0}\to \text{End}_S(A)$
such that the induced homomorphism
$O_{E_0}\to \text{End}_{O_S}(\text{Lie} A)=O_S$
is the same as that defined by the structure morphism
$S\to \text{Spec }E_0$.

(2)
$\bar k$ is an $H$-equivalent class of
an $\hat O_{E_0}$-isomorphism $k:T(A)\to T$
such that there exists a $ \hat\Bbb Z$-isomorphism $k'$
making the diagram 
$$\CD
\hat T(A)\times \hat T(A)@>>>
\hat\Bbb Z(1)\\
@V{k\times k}VV @VV{k'}V\\
\hat R\times \hat R
@>>{(x,y)\mapsto \text{Tr}_{E_0/\Bbb Q}(ax\bar y)}>
\hat\Bbb Z
\endCD$$
commutative.

\noindent
It is easily checked that the functor $N_{0,H,\hat R}$ is 
independent of a choice of $R$
upto uniquely determined canonical isomorphism.

By the theory of complex multiplication,
for a sufficiently small $H$,
the functor $N_{0,H,\hat R}$ is represented
by $N_H=\text{Spec }E_{0,H}$
where $E_{0,H}$ is the abelian extension
corresponding to the open subgroup
$H\subset \Bbb A_{E_0,f}^\times$
by the isomorphism 
$\Bbb A_{E_0,f}^\times/E^\times\simeq
\text{Gal}(E_0^{ab}/E_0)$
of class field theory.
Similarly as above,
a natural action of
$T_0(\Bbb A_f)=\Bbb A_{E_0,f}^\times$
on the projective systems
$N=(N_K)_K$
and on the universal CM elliptic curve
$b:A_0=(A_{0,\hat T,K})_{\hat T,K}\to N$
is defined.

\medskip

\subheading{6.1 Geometric construction on $M',N_0$}

We show that
the direct image
$R^1a_*\Bbb Q_\ell$ of
the universal abelian scheme
$a:A\to M'$ gives the
sheaf $\Cal F'$.
Using it, we construct the sheaf
$\Cal F^{\prime (k)}$ on $M'$
in a purely geometric way.
We will also define geometrically
$\Cal F(\chi)$
on $N_0$.

Let $K'\subset\hat O_D^\times$, $\hat T$ 
and  the universal $O_D$-abelian scheme
$a_{K'}:A_{K',\hat T}\to M'_{K'}$ 
be as in the modular interpretation 
in Section 5.
By the ring homomorphism
$O_D\to \text{End}_{M'_{K'}}(A_{K',\hat T})$,
we regard the direct image  $R^1a_{K'*}\Bbb Q_\ell$
as a sheaf of $D\otimes \Bbb Q_\ell$-modules
for every $\ell$.
It is independent of the choice of lattice $\hat T$.
A canonical action of
$G'(\Bbb A_f)$ is defined on the
system of sheaves $R^1a_{*}\Bbb Q_\ell=
(R^1a_{K'*}\Bbb Q_\ell)_{K'}$.
By the modular interpretation, it is easy to
see that the sheaf $R^1a_*\Bbb Q_\ell$
is isomorphic to the sheaf $\Cal F'$
with the action of $G'(\Bbb A_f)$
defined at the end of Section 4.
We will identify them in the following.

For each $i\in I$,
let $e_i\in D\otimes_{\Bbb Q}L$ be the
idempotent
defined at the end of Section 4.
We regarded $R^1a_* L_\lambda$
as a sheaf of $D\otimes_{\Bbb Q}L$-modules.
Then $e_i\in D\otimes_{\Bbb Q}L$ acts 
on it as a projector and
the $e_i$-part
$e_i\cdot R^1a_* L_\lambda$ 
is isomorphic to $\Cal F'_i$.
Since $D$ is generated by $1+pO_D$,
we may write each $e_i$
as an $L$-linear combination of
elements in $1+pO_D$.
Therefore $e_i$ is an
$L$-linear combination of
endomorphisms of $A$ over $M'$
whose degrees are prime to $p$.

One finds easily an idempotent
$e^{(k_i)}\in \Bbb Q[\Cal S_{w-2}]$
of the group algebra of a symmetric group such that
the $e^{(k_i)}$-part
$e^{(k_i)}\cdot {\Cal F'_i}^{\otimes w-2}$
is equal to $\text{Sym}^{k_i-2}\Cal F'_i
\otimes (\det \Cal F'_i)^{\otimes \frac {w-k_i}2}$.
The action of the symmetric group
$\Cal S_{w-2}$ on ${\Cal F'_i}^{\otimes w-2}$
is induced by
the action of it
on the fiber product $a^{w-2}:A^{w-2}\to M'$
over $M'$ as permutations.
One can also find easily
a $\Bbb Q$-linear combination $e^1$
of the multiplications by prime-to-$p$ integers
on $A$ such that
$e^1R^1a_*\Bbb Q_\ell=R^1a_*\Bbb Q_\ell$
and 
$e^1R^qa_*\Bbb Q_\ell=0$ for $q\neq 1$.

Taking their product,
we obtain
an algebraic correspondence $e'$
on the $(w-2)g$-fold self-fiber product
$A^{(w-2)g}$ of $A\to M'$ with coefficients in $L$
satisfying the following conditions.

(1) It is an $L$-linear
combinations of
permutations in $\Cal S_{g(w-2)}$
and 
endomorphisms of $A^{(w-2)g}$
as an abelian scheme over $M'$
whose degrees are prime to $p$.

(2) It acts as an idempotent on
the cohomology sheaf
$R^qa^{(w-2)g}_* L_\lambda$
where $a^{(w-2)g}$ denotes the map
$A^{(w-2)g}\to M'$.
We have
$e'R^qa^{(w-2)g}_* L_\lambda=
\Cal F^{\prime (k)}$
for $q=(w-2)g$ and 
$e'R^qa^{(w-2)g}_* L_\lambda=0$ otherwise.

Similarly, we construct $\Cal F(\chi)$.
Let $H\subset\hat O_{E_0}^\times$, $\hat R$ 
and  the universal $O_{E_0}$-elliptic curve
$b_H:A_{0,H,\hat R}\to N_{0,H}$ 
be as in the modular interpretation 
in Section 5.
By the ring homomorphism
$O_{E_0}\to \text{End}_{N_{0,H}}(A_{0,H,\hat R})$,
we regard the direct image  $R^1b_{H*}\Bbb Q_\ell$
as a sheaf of $E_0\otimes \Bbb Q_\ell$-modules
for every $\ell$.
It is independent of the choice of lattice $\hat R$.
A canonical action of
$\Bbb A_{E_0,f})$ is defined on the
system of sheaves $R^1b_*\Bbb Q_\ell=
(R^1b_{H*}\Bbb Q_\ell)_{H}$.
By the modular interpretation, it is easy to
see that the sheaf $R^1b_*\Bbb Q_\ell$
is isomorphic to the sheaf on $N_0$
associated to the inverse of the
tautological representation 
$E^\times\to GL_{\Bbb Q}(E): 
t\mapsto t^{-1}\times $.
We will identify them in the following.

Let $e_0\in E_0\otimes_{\Bbb Q}L$ be the
idempotent
corresponding to the inclusion $E_0\to L$.
We regarded $R^1a_{0,*} L_\lambda$
as a sheaf of $E_0\otimes_{\Bbb Q}L$-modules.
Then $e_0\in E_0\otimes_{\Bbb Q}L$ acts 
on it as a projector and
the $e_0$-part
$e_0\cdot R^1a_{0,*} L_\lambda$ 
is isomorphic to $\Cal F(\chi)$.
Similarly as above,
we may write each $e_0$
as an $L$-linear combination of
elements in $1+pO_D$.
Therefore $e_0$ is an
$L$-linear combination of
endomorphisms of $A_0$ over $N_0$
whose degrees are prime to $p$.
Similarly as above,
after modifing $e_0$ if necessary,
we also have
$e_0\cdot R^qb_{*} L_\lambda=0$
for $q\neq 1$.

\medskip
\subheading{6.2 Geometric construction on $M''$}

We extend the geometric construction on
$M'$ to $M''$.
We first study the relation between them.
Recall that
$G''=B^\times\times_{F^\times} E^\times$ and 
$G'$ is the inverse image of $\Bbb Q^\times\subset F^\times$
by $\nu=\text{Nrd}_{B/F}\times N_{E/F}: 
G^{\prime\prime} \to F^\times$.
For an open compact subgroup
$K''\subset G''(\Bbb A_f)$
and for $g\in G''(\Bbb A_f)$,
we put $K^{\prime g}=G'(\Bbb A_f)\cap gK^{\prime\prime} g^{-1}$.
Recall that
$g:M'_{K^{\prime g}}\to M''_{K''}$
denotes the composition
$M'_{K^{\prime g}}\to M''_{gK''g^{-1}}
@>g>> M''_{K''}$.
The double coset
$\tilde G\backslash G''(\Bbb A_f)/K''_1
=F^{\times +}\backslash
\Bbb A_{F,f}^\times/\Bbb A_{f}^\times\nu (K''_1)$
is finite.
If $\Sigma\subset G''(\Bbb A_f)$
is a complete set of representatives,
we have a finite etale surjection
$\amalg\ g:\coprod_{g\in \Sigma}
M'_{K^{\prime g}}\to M''_{K''}$.

\proclaim{Lemma 2}
Let $K''\subset G''(\Bbb A_f)$ be a compact open subgroup
and put $K'=K''\cap G'(\Bbb A_f)$.
Then for a sufficiently small open subgroup
$K''_1\subset K''$ containing $K'$
and for
a complete set $\Sigma$ of representatives of
the finite set
$\tilde G\backslash G''(\Bbb A_f)/K''_1$,
the map
$$\amalg\ g:\coprod_{g\in \Sigma}
M'_{K^{\prime g}}\to M''_{K''_1}$$
is an isomorphism.
\endproclaim
\demo{Proof}
Since it is an etale surjection,
it is enough to show the map is injective
on the $\Bbb C$-valued points.
Since $\Sigma$ is a complete
set of representatives,
it is enough to consider each map $g$.
Let $\bar \nu:G''(\Bbb A_f)\to \Bbb A_{F,f}^\times/
\Bbb A_{\Bbb Q,f}^\times$
denote the map induced by $\nu$.
We claim that the equality
$\bar \nu(K'')\cap (\overline{O_F^\times}/\Bbb Z^\times)
=\bar \nu(K''\cap \overline{O_E^\times})$
implies the injectivity of the map
$g:M'_{K^{\prime g}}(\Bbb C)\to M''_{K''}(\Bbb C)$.

We prove Lemma admitting the claim.
Namely, we prove that for a sufficiently
small open subgroup $K''_1\supset K'$ of $K''$,
we have an equality
$\bar \nu(K''_1)\cap (\overline{O_F^\times}/\Bbb Z^\times)
=\bar \nu(K''_1\cap \overline{O_E^\times})$.
Since $N_{E/F}(O_E^\times)$ is of finite index in
$O_F^\times$, the right hand side
$\bar \nu(K''\cap \overline{O_E^\times})$
is an open subgroup of
the left hand side
$\bar \nu(K'')\cap (\overline{O_F^\times}/\Bbb Z^\times)$.
Hence,
for a sufficiently small open subgroup $\overline{K_1}$
of $K''/K'\simeq \bar \nu (K'')$ we have
$\overline{K_1}\cap
 (\overline{O_F^\times}/\Bbb Z^\times)
=\overline{K_1}\cap
\bar \nu(K''\cap \overline{O_E^\times})$.
For the corresponding open subgroup $K''_1=
K''\cap \bar\nu^{-1}(\overline{K_1})$,
this is nothing but the required equality
$\bar \nu(K''_1)\cap (\overline{O_F^\times}/\Bbb Z^\times)
=\bar \nu(K''_1\cap \overline{O_E^\times})$.

We prove the claim.
Namely, we assume
$\bar \nu(K'')\cap (\overline{O_F^\times}/\Bbb Z^\times)
=\bar \nu(K''\cap \overline{O_E^\times})$
and prove the map
$g:M'_{K^{\prime g}}(\Bbb C)\to M''_{K''}(\Bbb C)$
is injective.
Replacing $K''$ by $gK''g^{-1}$,
it is enough to show that the map
$M'_{K'}(\Bbb C)\to M''_{K''}(\Bbb C)$
is injective
for $K'=K''\cap G'(\Bbb A_f)
=\text{Ker}(\bar \nu:K''\to 
\hat O_F^\times/\hat\Bbb Z^\times)$.
We consider the commutative diagram of exact sequences
$$\CD
@.K''\cap \overline{O_E^\times}
@>>> K''/K' @>>>
K''/(K''\cap \overline{O_E^\times})K'@>>>1\\
@.@V{\bar\nu}VV @V{\bar\nu}V{\cap}V @VVV @.\\
1@>>>\overline{O_F^\times}/\Bbb Z^\times
@>>>\hat O_F^\times/\hat \Bbb Z^\times
@>>>\hat O_F^\times/\hat \Bbb Z^\times
\overline{O_F^\times}.
\endCD$$
The middle vertical arrow is injective by
the definition of $K'$.
By the snake lemma,
the equality
$\bar \nu(K'')\cap (\overline{O_F^\times}/\Bbb Z^\times)
=\bar \nu(K''\cap \overline{O_E^\times})$
is equivalent to
the injectivity of the right vertical arrow.
Since $\hat O_F^\times/\hat \Bbb Z^\times\overline{O_F^\times}$
is a subgroup of 
$\Bbb A_{F,f}^\times/\Bbb A_{\Bbb Q,f}^\times\overline{F^\times}$,
we get an exact sequence
$$K'/(K'\cap \overline{O_E^\times})@>>>
K''/(K''\cap \overline{O_E^\times})@>>>
\Bbb A_{F,f}^\times/\Bbb A_{\Bbb Q,f}^\times\overline{F^\times}.$$
We consider the commutative
diagram
$$\CD
M'(\Bbb C)=\varprojlim_{K'}M'_{K'}(\Bbb C)@>>>
M''(\Bbb C)=\varprojlim_{K''}M''_{K''}(\Bbb C)\\
@VVV @VVV\\
\Bbb Q^\times\backslash \Bbb A_{\Bbb Q,f}^\times
@>>>
\overline{F^\times}\backslash \Bbb A_{F,f}^\times.
\endCD$$
The horizontal arrows are injective
by Variante 1.15.1 and Lemma 1.15.3 [D2].
We have 
$M'_{K'}(\Bbb C)=M'(\Bbb C)/K'$ and
$M''_{K''}(\Bbb C)=M''(\Bbb C)/K''$.
From these facts, it is straightforward to
show that the exactness above implies
the injectivity of the canonical map 
$M'_{K'}(\Bbb C)\to
M''_{K''}(\Bbb C)$.
\enddemo

We extend the universal $O_D$-abelian scheme
$A$ on $M'$ to an $O_D$-abelian scheme also denoted
by $A$ on $M''$.
Let $K''\subset G''(\Bbb A_f)$ be a sufficiently small
open subgroup.
We assume that the map 
$\amalg\ g:\coprod_{g\in \Sigma}
M'_{K^{\prime g}}\to M''_{K''_1}$
in Lemma 2 
is an isomorphism.
We take an $\hat O_D$-lattice $\hat T$
in $D\otimes\Bbb A_f$.
For each $g\in \Sigma$,
we have a $gO_Dg^{-1}$-abelian scheme $A_{K'_g, g\hat T}$
on $M'_{K'_g}$ since $g\hat T$ is $K'_g$-stable.
We define an abelian scheme $A_{K'',\hat T}$
on $M''_{K''}$ to be $A_{K'_g,g\hat T}$
on the image of $M'_{K'_g}$.
We define an $O_D$-multiplication
on $A_{K'',\hat T}$ as
$O_D@>{a\mapsto gag^{-1}}>> gO_Dg^{-1}
\to \text{End}_{M'_{K'_g}}(A_{K'_g,g\hat T})$
on $M'_{K'_g}$.
By the action of $\tilde G$
described in the previous section,
we see that the abelian scheme  $A_{K'_g, g\hat T}$
is independent of the choice of representatives
$\Sigma$.
We also see
by the action of $\tilde G$ that,
for $g\in G''(\Bbb A_f)$,
compact open subgroups
$K''_1,K''_2\subset G''(\Bbb A_f)$
and $K''_i$-stable $\hat O_D$-lattices $\hat T_i$
satisfying $g^{-1}K_1''g\subset K''_2$,
we have an isogeny
$A_{K''_1,\hat T_1}\to g^* A_{K''_2,\hat T_2}$.
Thus we obtain an action of $G''(\Bbb A_f)$
on the projective system
$A=(A_{K'',\hat T})_{K'',\hat T}$
over $M''=(M''_{K''})_{K''}$.

On the $(w-2)g$-fold self-fiber product
$A^{(w-2)g}$ of $A\to M''$,
we define 
an algebraic correspondence $e'$
with coefficients in $L$
exactly in the same way as in the case of $M'$.
Then, it is an $L$-linear
combinations of
permutations in $\Cal S_{g(w-2)}$
and 
endomorphisms of $A^{(w-2)g}$
as an abelian scheme over $M'$
whose degrees are prime to $p$.
Further,
it acts as an idempotent on
the cohomology sheaf
$R^qa^{(w-2)g}_* L_\lambda$
where $a^{(w-2)g}$ denotes the map
$A^{(w-2)g}\to M'$.
We have
$$e'R^qa^{(w-2)g}_*L_\lambda=
\bigotimes_i
\text{Sym}^{k_i-2}(e_i\cdot R^1a_* L_\lambda)
\otimes (\det e_i\cdot R^1a_* L_\lambda
)^{\otimes \frac {w-k_i}2}
$$ 
for $q=(w-2)g$ and 
$e'R^qa^{(w-2)g}_*L_\lambda=0$ otherwise.
By the modular interpretation of
$M'$, we see that the $K''$-equivalent class of
the isomorphism $\hat T\to T(A_{K'',\hat T})$ is well-defined.
Passing to the limit,
we obtain an isomorphism
$D\otimes \Bbb A_f\to 
R^1a_*\Bbb Q_\ell$ on
$\varprojlim_{K''}M_{K''}$.
The isomorphism is
compatible with the action of
$G''(\Bbb A_f)$.
On the left hand side
$D\otimes \Bbb A_f$,
the group $G''(\Bbb A_f)\subset 
(D\otimes \Bbb A_f)^\times$
acts by the multiplication by the inverse of
the main involution:
$t\mapsto \bar t^{-1}\times$.
Thus similarly as on $M'$,
we have
$$e'R^{(w-2)g}a^{(w-2)g}_* L_\lambda=
\Cal F^{\prime \prime (k)}.$$

\medskip
\subheading{6.3. Geometric construction on $M$}

We will define an analogue 
$c:X\to M\times N$ of Kuga-Sato variety
and an algebraic correspondence
$e=e^{(k)}$ on $X$
with coefficient in $L$
satisfying the following property:
It is an $L$-linear combination
of endomorphisms of $X$ as
an abelian scheme over $M\times N$
whose degree are prime to $p$.
It acts as an idempotent on
the higher direct image 
$R^qc_*\Bbb Q_\ell\otimes L
=\prod_{\lambda|\ell}R^qc_*L_\lambda$.
The image of the projector
$e\cdot R^qc_*L_\lambda$
is a smooth $L_\lambda$-sheaf
isomorphic to 
$$\alpha^*\Cal F^{(k)}
\otimes \beta^* F(\chi)^{\otimes(w-2)(g-1)}=
pr_1^*\Cal F^{(k)}
\otimes \beta^* F(\chi_0)^{\otimes(w-2)(g-1)}$$
for $q=q_0=(2g-1)(w-2)$
and is 0 otherwise.

We define $X$ to be the fiber product
$$X=\alpha^*A^{g(w-2)}\times_{M_E\times N}
\beta^*A_0^{(g-1)(w-2)}.$$
Here $\alpha^*A^{g(w-2)}$ denotes the base change by
$\alpha:M\times N\to M''$ of the
$g(w-2)$-fold self fiber product of $A\to M''$.
Similarly
$\beta^*A_0^{g(w-2)}$
denotes the base change by
$\beta:M\times N\to N_0$ of the
$(g-1)(w-2)$-fold self fiber product of $A_0 \to N_0$.
The symbol $X$ denotes
the projective system
$X=(X_{K,H,\hat T,\hat R})_{K,H,\hat T,\hat R}$
of abelian schemes over 
$M\times N=(M_K\times N_{H})_{K,H}.$

Next 
we define an algebraic correspondence 
$e=e^{(k)}$ on $X$.
We have defined algebraic correspondences
$e'$ on $A^{g(w-2)}$ on $M''$
and $e_0$
on $A^{g(w-2)}$ on $A_0$ on $N_0$
and the end of subsections 6.2 
and 6.1 respectively.
Let $e_0^{\otimes (g-1)(w-2)}
=\prod_{i=1}^{(g-1)(w-2)}pr_ie_0$
be the algebraic correspondence
on the $(g-1)(w-2)$-nd self fiber product
$A_0^{(g-1)(w-2)}$
defined as the product of the pull-back
of the algebraic correspondence $e_0$ on $A_0$
by projections.
We define an algebraic correspondece $e$
on $X$ as the product of the pull-back of $e'$
by $\alpha$
with the pull-back of $e_0^{\otimes (g-1)(w-2)}$
by $\beta$.
Namely we put
$e=\alpha^*e'\times
\beta e_0^{\otimes (g-1)(w-2)}.$
Then it satisfies the required property
stated in the beginning of this section.

Let $H\subset \Bbb A_{E,f}^\times$ be
a sufficiently small open compact subgroup.
Let $\frak m=\frak n O_E$ be a
sufficiently divisible integral ideal of $O_E$.
We assume $H=
H^{\frak m}\cdot H_{\frak m}$
is the product of the
prime-to-$\frak m$ component
$H^{\frak m}=\prod_{\frak s\nmid \frak m}
O_{E,\frak s}^\times$
with the $\frak m$-primary component
$H_{\frak m}$.
Let
$T_0^{\frak m}=
L[P_{\frak s};\ {\frak s}\nmid {\frak m}]$
be the free $L$-algebra
generated by the class $P_{\frak s}$
of the
inverse of prime element for
${\frak s}\nmid {\frak m}$.
We consider
$H^q(X_{K,H,\hat T,\hat R}\otimes_E \bar E,L_\lambda)$
as a $T^{\frak n}
\times T_0^{\frak m}$-module
and 
$H^0(N_{H,\bar E},\Cal F(\chi_0))$
as a $T_0^{\frak m}$-module.

Applying the Leray spectral sequence
to $c:X_{K,H,\hat T,\hat R}
\to M_K\times_FN_{H}$,
we obtain the following.

\proclaim{Lemma 3}
Let 
$K\subset G(\Bbb A_f)$
and $H\subset \Bbb A_{E,f}^\times$
be sufficiently small open compact subgroups
and let $\hat T\subset V\otimes \Bbb A_f$ and
$\hat R \subset E_0\otimes \Bbb A_f$
be an $\hat O_D$-lattice and
an $\hat O_{E_0}$-lattice respectively.
Let $X=X_{K,H,\hat T,\hat R}$
be the analogue of Kuga-Sato variety defined above.
Then there is an algebraic correspondence
$e$ on $X$ with coefficient
in $L$ satisfying the following
properties.

\noindent (1)
There exists
elements  $a_i\in L$,
permutations $\tau_i
\in \Cal S_{g(w-2)}$ of
the first $g(w-2)$-factors in $X$
and
endomorphisms $\varphi_i\in \text{End}_MX$
of degree prime to $p$
such that
$$e=\sum_i a_i\tau_i\varphi_i.$$

\noindent (2)
For each finite place $\lambda$ of $L$,
the action of $e$ on 
$H^q(X_{K,H,\hat T,\hat R}
\otimes_E\bar E,L_\lambda)$
is a projector.
Put $q_0=(2g-1)(w-2)$.
Then, there is an isomorphism
$$\align
&\ e\cdot H^q(X_{K,H,\hat T,\hat R}\otimes_E\bar E,L_\lambda)\\
\simeq&\ 
H^{q-q_0}(M_K\otimes_F\bar F,\Cal F_\lambda^{(k)})
\otimes_{\Bbb Q_\ell}
H^0(N_{H}\otimes_E\bar E,\Cal F(\chi_0^{(g-1)(w-2)})).
\endalign
$$
The isomorphism is compatible with
the actions of the
absolute Galois group
$G_E=\text{Gal}(\bar E/E)$
and of the Hecke algebra
$T^{\frak n}\otimes T_0^{\frak m}$.
\endproclaim

Using Lemma 3,
we give a statement, Claim 3, in terms of $X$
and $e$,
implying Claim 2 and hence Theorems.
Let $\frak q$ be a place of $E$
dividing $\frak p$.
By the assumption that $p$
splits in $E_0$,
the local field 
$E_{\frak q}$ is
canonically isomorphic to $F_{\frak p}$.
We identify 
$F_{\frak p}=E_{\frak q}$
by the canonical isomorphism.
Since we want to prove the assertions
on the action of Galois group
$\text{Gal}(\bar F_{\frak p}/F_{\frak p})$,
it is enough to consider the action of
$\text{Gal}(\bar E_{\frak q}/E_{\frak q})$,
induced by the isomorphism.

\proclaim{Claim 3}
We keep the notation in Claim 2.
Let $K\subset G(\Bbb A_f)$
and
$H\subset \Bbb A_{E,f}^\times$
be sufficiently small
open compact subgroups.
Let $X=X_{K,H,\hat T,\hat R}$
denote the analogue of Kuga-Sato variety.
Then, the following holds.

\noindent (0)
The $p$-adic representation
$H^q(X\otimes_E\bar E_{\frak p},\Bbb Q_p)$
of $G_{E_{\frak p}}=
\text{\rm Gal}(\bar E_{\frak p}/E_{\frak p})$
is potentially semi-stable for all $q$.

\noindent (1)
Let $\sigma\in W^+=\{\sigma\in W(\bar E_{\frak p}/E_{\frak p})|
n(\sigma)\ge0\},
T\in T^{\frak n},
P \in T_0^{\frak m},
\tau\in \Cal S_{g(w-2)}$
and let $\psi:X\to X$
be an endomorphism of degree prime to $p$.
Then for the composite
$\Gamma =T\circ R \circ
\tau\circ \psi$ as an algebraic correspondence,
we have an equality in $\Bbb Q$
$$\sum_q(-1)^q\text{Tr}(\sigma\circ \Gamma|
H^q(X\otimes_E\bar E_{\frak p},\Bbb Q_\ell))
=\sum_q(-1)^q\text{Tr}(\sigma\circ \Gamma|
D(H^q(X\otimes_E\bar E_{\frak p},\Bbb Q_p))).$$

\noindent (2)
Let $e$ be the algebraic correspondence in Lemma 3
and let $\mu|p$ be a finite place of $L\supset E_0$.
Then the monodromy filtration of
the representations
$e\cdot
H^q(X\otimes_E\bar E_{\frak q},L_\lambda)$
and
$D(e\cdot
H^q(X\otimes_E\bar E_{\frak q},L_\mu))$
of the Weil-Deligne group
${}'W(\bar E_{\frak q}/E_{\frak q})$
are pure of weight $q$.
\endproclaim

We deduce each assertion in Claim 2 from 
the corresponding assertion in Claim 3.
Since we identify $F_{\frak p}=
E_{\frak q}$,
it is sufficient to consider the representations
of the Weil-Deligne group
${}'W(\bar E_{\frak q}/E_{\frak q})$.
By Lemma 3, the representation
$H^q(M_{\bar E_{\frak q}},\Cal F_\lambda^{(k)})$
is a direct summand of
$e\cdot H^{q+q_0}(X_{\bar E_{\frak p}},L_\lambda)
((g-1)(w-2))$.
Therefore the assertions (0) and (2) in Claim 2
follows from the assertions (0) and (2) in Claim 3
respectively.
We deduce the assertion (1) in Claim 2 from
the assertion (1) in Claim 3.
By the description of
$\Cal F(\chi_0)$
given at the end of Section 4,
we find easily an element $e^\circ
\in 
T_0^{\frak m}$
acting as a projector
$H^0(N_{H}\otimes_E\bar E,
\Cal F(\chi_0^{(g-1)(w-2)}))
\to L_\lambda(-(g-1)(w-2))$.
Thus by Lemma 3, 
there is an isomorphism
$$
e^\circ\circ
e\cdot H^q(X_{K,K^{\circ}}\otimes_E\bar E,L_\lambda)
\to
H^{q-q_0}(M_K\otimes_F\bar F,\Cal F_\lambda^{(k)})
(-(g-1)(w-2))
$$
compatible with the actions of
the Galois group 
$G_E=\text{Gal}(\bar E/E)$
and of the Hecke algebra
$T^{\frak n}$.
Hence the equality in (1) in Claim 3
implies the equality in (1) in Claim 2.
Thus Theorems 1 and 2,
are reduced to Claim 3.

We may deduce the assertion (0)
using alteration [dJ].
We will give a proof without using alteration
by constructing a semistable model
of $X$.

For later use, we describe the
Hecke operaters $T_{\frak r}\in T^{\frak n}$
and $P_{\frak s}\in T_0^{\frak m}$
for primes $\frak r\nmid \frak n$
of $O_F$ and
$\frak s\nmid \frak m$
of $O_E$ respectively.
Write
$X=X_{K,H,\hat T,\hat T^\circ}$
and 
$M\times N=M_K\times N_{H}$
for short.
For $\frak r$,
it is defined 
as $T_{\frak r}=
p_{1*}\circ q^* \circ p_2^*$
where $p_1,p_2,q$ are as in 
the diagram
$$\CD
X
@<{p_1}<<
X_{K_g,H,\hat T,\hat R}
@>q>>
X_{K_g,H,g\hat T,\hat R}
@>{p_2}>>
X\\
@VVV @VVV @VVV @VVV\\
M\times N
@<{p_1}<<
M_{K_g}\times N_{H}
@=
M_{K_g}\times N_{H}
@>{p_2}>>
M\times N.
\endCD$$
In the diagram,
$g=g_{\frak r}\in G(\Bbb A_f)$
is an element whose $\frak r$-component
is $\pmatrix \pi_{\frak r}^{-1} &0 \\ 0& 1\endpmatrix$
and other components are 1
and
$K_g=K\cap g Kg^{-1}$.
The map $p_1$ is induced by the
inclusion $K_g\to K$,
the map $p_2=g_*$ is induced by $g$
and the left and right squares are cartesian.
The map
$q$ is an isogeny corresponding to
the inclusion $\hat T\to g\hat T$.

Similarly for $\frak s$,
the operater is defined 
as $P_{\frak s}=
q^* \circ p_2^*$
where $p_2,q$ are as in 
the diagram
$$\CD
X
@>q>>
X_{K,H,g\hat T,N_{E/E_0}g\hat R}
@>{p_2}>>
X\\
@VVV @VVV @VVV\\
M\times N
@=
M\times N
@>{p_2}>>
M\times N.
\endCD$$
In the diagram,
$g\in \Bbb A_{E,f}^\times$
denotes an element whose $\frak s$-component
is the inverse of a prime element
$\pi_{\frak s}^{-1}$
and the other components are 1.
The map $p_2=g_*$ is induced by $g$
and right square is cartesian.
The map 
$q$ is an isogeny corresponding to
the inclusions $\hat T\to g\hat T$ and
$\hat R \to N_{E/E_0}g\cdot \hat R$.

\medskip

\subheading{7. Semi-stable model}

In the last section,
we defined an analogue of Kuga-Sato
variety as an abelian schemes on 
a Shimura curve.
The goal, Lemma 4, in this section is
to give their semi-stable models.

We prepare some terminology.
Let $K,H,
\hat T,\hat R$
be as in the last section.
We assume that
each component of the generic fiber
$M_K\otimes_F\bar F$ are of genus greater than 1.
Then by the stable reduction theorem of curve,
for a sufficiently large finite extension $V$
of the maximal unramified extension
$\widehat{F_{\frak p}^{nr}}$,
the base change
$M_{K,V}=M_K\otimes_FV$
admits a semi-stable model (not necessarily connected)
over the integer ring $O_V$.
We do not need unramified extension but
for later use, we will do it here.
We take the minimal one among the semi-stable models
over $O_V$ and denote it by
$M_{K,O_V}$.
Recall that we identified the local field
$E_{\frak p}$ with
$F_{\frak p}$.
From now on,
we consider $V$ as
an extension of 
$E_{\frak p}$ by this identification.
Since $N_{H}$ is the disjoint union
of the spectrum of finite extensions of $E$,
the base change
$(M_K\times_F N_{H})\otimes_EV$
also admits a semi-stable model over the
integer ring $O_V$.
We also take the minimal one among them and
write it by 
$(M_K\times_F N_{H})_{O_V}$.
We claim the following.

\proclaim{Lemma 4}
Let $K,H$ and $V$ be as above.

(1)
Let $g\in G(\Bbb A_f), h\in \Bbb A_{E,f}^\times$
and let $K_1\subset gKg^{-1},H_1\subset H$
be open compact subgroups.
Assume that 
the groups are of the form
$K=K_{\frak p}
K^{\frak p},
g^{-1}K_1g=K_{\frak p}
(g^{-1}K_1g)^{\frak p},
H=H_{\frak q}
H^{\frak q}$
and 
$H_1=H_{\frak q}
H_1^{\frak q}$.
Then the pull-back
of the map
$(g,h)_*:M_{K_1}\times_FN_{H_1}\to
M_K\times_F N_H$ to
the base change over $V$
is extended uniquely to a finite etale morphism
$(g,h)_*:(M_{K_1}\times_FN_{H_1})_{O_V}\to
(M_K\times_F N_H)_{O_V}$
of the minimal semi-stable models.

(2)
Let $\hat T,\hat R$ be as above.
Then the pull-back
of the abelian scheme
$X_{K,H,\hat T,\hat R}
\to
M_K\times_F N_{H}$ to
the base extension
$(M_K\times_F N_{H})\otimes_EV$
is extended uniquely to an abelian scheme
over a semi-stable model.

(3)
Let $\hat T_1,\hat R_1$
be sublattices in
$\hat T$ and $\hat R$ in (2)
respectively. Assume that
their $p$-components are the same.
Then the pull-back
of the isogeny
$X_{K,H,\hat T_1,\hat R_1}
\to
X_{K,H,\hat T,\hat R_1}$
on $M_K\times_F N_{H}$ to
the base extension
$(M_K\times_F N_{H})\otimes_EV$
is extended uniquely to an etale isogeny
over a semi-stable model.
\endproclaim

\demo{Proof}
(1). We may assume $g=1$ and $h=1$.
Further we may assume $H=H_1$.
In fact, the map $N_{H_1}\to N_H$ is
unramified at $\frak q$
by the assumption that
there $\frak q$-components are the same
by class field theory.
Further, it is sufficient to
show that the map
$M_{K_1}\to M_K$ is extended to
a finite etale morphism of
minimal semi-stable models
$M_{K_1,O_V}\to M_{K,O_V}$.
In fact,
then the fiber product
$M_{K_1,O_V}
\times_{ M_{K,O_V}}
(M_{K}\times N_H)_{O_V}$
is a semi-stable model of
$(M_{K_1}\times N_H)_V$
and does not have a $(-1)$-curve.
Hence it is the minimal semi-stable model
$(M_{K_1}\times N_H)_{O_V}$
and 
$(M_{K_1}\times N_H)_{O_V}\to
(M_{K}\times N_H)_{O_V}$
is finite etale.

In the case where
the $\frak p$-components of
$K_{\frak p}=K_{1,\frak p}$
are $GL_2(O_{F,\frak p})$,
it is shown in [C1] Proposition 6.1 
and in loc.\ cit.\ 6.2 that
the canonical map
$M_{K_1,F_{\frak p}} \to M_{K,F_{\frak p}}$
is extended to a finite etale morphism
$M_{K_1,O_{F,\frak p}}\to M_{K,O_{F,\frak p}}$
of proper smooth models.
We consider the general case.
Let $\bar K\supset K,\bar K_1\supset K_1$
be the groups obtained by replacing
their $\frak p$-components
$K_{\frak p}=K_{1,\frak p}$
by $GL_2(O_{F,\frak p})$
respectively.
First, we show that
the canonical map 
$M_{K,V}
\to M_{\bar K,V}$
is extended to
the minimal semi-stable model
$M_{K,O_V}
\to M_{\bar K,O_V}$.
In fact, it is extended on a suitable blow-up.
However, the exceptional divisors are contracted
to points in the image and hence the
map is defined on the semi-stable model.
We consider the fiber product
$M_{\bar K_1,O_V}
\times_{ M_{\bar K,O_V}}
M_{K,O_V}$.
It is a semi-stable model of
$M_{K_1,V}$
and does not have a $(-1)$-curve.
Hence it is minimal and the assertion is proved.

(2).
We assume there
exists an open compact subgroup
$K''\subset G''(\Bbb A_f)$ 
containing $K''\supset KH$ 
and satisfying
the following conditions (a) and (b).

(a) $K''$ satisfies the conclusion
of Lemma 2.
Namely for a complete set $\Sigma$ of representatives
$\tilde G\backslash G''(\Bbb A_f)/K''$,
the map $\amalg g:\coprod_gM'_{K'_g}\to M''_{K''}$
is an isomorphism.

\noindent
To state the other condition,
we identify the group
$G'(\Bbb Q_p)$.
By the assumption that $E_0$ splits at $p$,
we have an isomorphism
$$\CD
G'(\Bbb Q_p)@=
\Bbb Q_p^\times\times (B\otimes_{\Bbb Q}\Bbb Q_p)^\times
@>\sim>>\Bbb Q_p^\times\times GL_2(F_{\frak p})\times
(B\otimes_{F}F_p^{\frak p})^\times\\
@V{\cap}VV @V{\cap}VV @.\\
G''(\Bbb Q_p)
@=
(F\otimes \Bbb Q_p)^\times\times (B\otimes_{\Bbb Q}\Bbb Q_p)^\times
@.\endCD$$
[C1] (2.6.3).
The second condition is the following.

(b) The intersection
$K'=K''\cap  G'(\Bbb A_f)$
is of the form
$K'=\Bbb Z_p^\times\times GL_2(O_{F,\frak p})
\times K_p^{\prime {\frak p}}\times K^{\prime p}$
for some choice of isomorphism as above.

\noindent
It is shown in [C1] Proposition 5.4
using a modular
interpretation that
the condition (b) implies that
$M'_{K'_g}$ has good reduction
over $O_{E,\frak p}$ and the abelian scheme
$A_{K'_g,g\hat T}$ on
the generic fiber is extended
to a (unique) proper smooth
model $M'_{K'_g,O_{E,\frak p}}$.
We will recall the modular interpretation
in Section 9.
Hence by the condition (a),
$M''_{K''}$ has also good reduction
over $O_{E,\frak p}$ and the abelian scheme
$A_{K'',\hat T}$ is extended to the 
proper smooth model $M''_{K'',O_{E,\frak p}}$.
By the same argument as in the proof of (1),
the map $(M_K\times N_{H})_V\to
M''_{K'',V}$ is extended 
uniquely to a map $(M_K\times N_{H})_O\to
M''_{K'',O_{E,\frak p}}\otimes O_V$.
Hence we obtain the extension of
an abelian scheme
$X_{K,H,\hat T,\hat R}$,
by taking the pull-back.

(3).
Since $(M\times N)_O$
is normal,
an endomorphism on the generic fiber
is extended to the integral model
by a theorem
of Grothendieck [G].
\enddemo
In the proof of (3),
we could also use the modular interpretation
recalled in Section 9.

\medskip
\subheading{8. Proof of Theorems}

We prove Theorems 1 and 2
by showing the assertions in Claim 3.
The argument is the same as
in [Sa].
Let the notation be as in Claim 3.
We fix sufficiently small open compact subgroups 
$K\subset G(\Bbb A_f),
H\subset \Bbb A_{E_0,f}^\times$,
an $\hat O_D$-lattice $\hat T$
and $\hat O_{E_0}$-lattice $\hat R$.
To simplify the notation,
we will write
$M\times N$ for $M_K\times N_H$ and
$X$ for $X_{K,H,\hat T,\hat R}$.
Recall that we identify
$E_{\frak q}=F_{\frak p}$.

We prove that
the $p$-adic representation
$H^q(X\otimes_E\bar E_{\frak q},\Bbb Q_p)$
of the Galois group 
$\text{Gal}(\bar E_{\frak q}/E_{\frak q})$
is potentially semi-stable.
Since we have a semi-stable
model $X_{O_V}$ of 
the base change $X_V$ to an extension
$V$ of $E_{\frak q}$ by Lemma 4,
we just apply
the C$_{st}$-conjecture proved by Tsuji [Tj]
to a semi-stable scheme $X_O$.

We compute
$D_{pst}(H^q(X\otimes_E\bar E_{\frak q},\Bbb Q_p))$
in terms of the minimal semi-stable model 
$X_O$ of $X$
defined in Lemma 4.
Let $Y$
denote the closed fiber of the 
minimal semi-stable model $X_O$
with the natural log structure.
Then further by [Tj],
we have a canonical isomorphism
$$D_{pst}(H^q(X\otimes_E\bar E_{\frak q},\Bbb Q_p))
\simeq
H^q_{\text{log crys}}(Y/W)\otimes_{\Bbb Z_p}
\Bbb Q_p.$$
It follows from the functoriality
of the comparison isomorphism for finite etale morphism
and from the compatibility with the Poincar\'e duality
that the isomorphism is compatible with
the action of endomorphisms
and permutations appeared in Claim 3.
We define Hecke operators
on the log cristalline cohomology
and compare them with those
on the left hand side
induced by the Hecke operators
on the etale cohomology.
Let $\frak n\subset O_F$
and $\frak m\subset O_E$
be sufficiently divisible ideals
as in Section 6.3.
Let $\frak r\nmid \frak n$ 
be a prime ideal of $O_F$.
Then the projections $p_1,p_2$ and the
isogeny $q$ described at the end of
section 7 is extended to a finite etale
morphism of the
minimal semi-stable model by Lemma 4.
On log cristalline cohomology,
we define the Hecke operator $T_{\frak r}$
as the composite $p_{1*}\circ q^*\circ p_2^*$.
Similarly we define
the Hecke operator $P_{\frak s}$
for a prime ideal $\frak s\nmid \frak m$ 
of $O_E$ as the composite $q^*\circ p_2^*$
using the description at the end of
section 7.
Then it follows from the functoriality
that the isomorphism is compatible with the
Hecke operators thus defined.

We define the Galois action
on the log cristalline cohomology
and compare it with that on
the left hand side defined in Section 2.
We may and do assume that the
finite extension $V$ of $\widehat{E_{\frak q}^{nr}}$
is the completion of a Galois extension of $E_{\frak q}$.
We have a natural action of the Galois group 
$G_{\frak q}=\text{Gal}(\bar E_{\frak q}/E_{\frak q})$
on $V$ and hence on the base change
$M_{V}$.
Since the minimal semi-stable model is unique,
the action of
$G_{\frak q}$
on the generic fiber $M_{V}$ is extended to the
minimal semi-stable model $M_{O_V}$.
Further it is uniquely extended to
that on the abelian scheme $X_{O_V}$.
It induces a semi-linear action of the Weil group
$W_{{\frak q}}$ on the log cristalline cohomology.
By modifying the action of 
$\sigma\in W_{E_{\frak q}}$
by
$\varphi^{n(\sigma)}\circ \sigma$
as in Section 2
and together with $N$,
we define a linear action of 
the Weil-Deligne group $W'_{E_{\frak q}}$
on the log cristalline cohomology.
We verify the compatibility of the isomorphism
with the action of Weil-Deligne group
defined above.
By transport of the structure,
it is compatible with the semi-linear
action of the Weil group
before modification.
Since the comparison isomorphism is
compatible with the action of
$F$ and $N$,
the compatibility is established.

Therefore, 
Claim 3 is reduced to the following.
\proclaim{Claim 4}
Let the notation be as in Claim 3.
Then,
the following holds.

(1)
Let $\sigma\in W^+=\{\sigma\in 
W(\bar E_{\frak q}/E_{\frak q})|
n(\sigma)\ge0\},
T\in T^{\frak n},
R\in T_0^{\frak m},
\tau \in \Cal S_{w-2}^g$
and let $\psi:X\to X$
be an endomorphism of degree prime to $p$.
Then for the composite
$\Gamma =T\circ R \circ
\tau\circ \psi$ as an algebraic correspondence,
we have an equality in $\Bbb Q$
$$\sum_q(-1)^q\text{Tr}(\sigma\circ \Gamma|
H^q(X\otimes_E\bar E_{\frak q},\Bbb Q_\ell))=
\sum_q(-1)^q\text{Tr}(\sigma\circ \Gamma|
H^q_{\text{log crys}}(Y/W)).$$

(2)
Let $e$ be the algebraic correspondence in Lemma 3
and let $\lambda\nmid p,\mu|p$
be finite places of $L\supset E_0$.
Then the monodromy filtration of
the representations
$e\cdot
H^q(X\otimes_E\bar E_{\frak q},L_\lambda)$
and
$e\cdot
(H^q_{\text{log crys}}(Y/W)\otimes
\widehat{L_\mu^{nr}})$
of the Weil-Deligne group
${}'W(\bar E_{\frak q}/E_{\frak q})$
are pure of weight $q$.
\endproclaim

In (2), the tensor product is taken with
respect to the map
$W=O_{\widehat{E_{\frak q,0}^{nr}}}
\subset \widehat{E_{\frak q,0}^{nr}}
@<\sim<< \widehat{F_{\frak p,0}^{nr}}
\to \widehat{L_{\mu}^{nr}}$
where the last map is that fixed in section 2.
We prove Claim 4
by studying the weight spectral sequences
for $\ell\neq p$ and for $p$.

We prove (1).
First we compute the $\ell$-adic case.
We consider the weight spectral sequence [RZ], [I]
$$E_1^{i,j}=
\bigoplus_{k\ge \max(0,-i)}
H^{j-2k}(Y^{(i+2k)},\Bbb Q_\ell(-k))\Rightarrow
H^{i+j}(X_{\bar V},\Bbb Q_\ell)$$
for the semi-stable model $X_{O_V}$.
Here $Y^{(i)}$ denotes the disjoint union of
$i+1$ by $i+1$ intersections of the
irreducible components of the closed fiber
$Y=X\otimes_{O_{E',{\frak q}}} \bar\Bbb F_{\frak q}$.
The schemes $Y^{(i)}$ are projective and smooth over
$\bar\Bbb F_{\frak q}$.
Since the action of the Galois group $G_{\frak q}$
is extended to the semi-stable model $X_{O_V}$,
the spectral sequence is compatible
with its action by transport of structure.
It is also compatible with 
the action of Hecke operators,
endomorphisms and permutations
by the same argument as
in the case of $p$-adic
comparison isomorphism above.
Hence the left hand side of the equality
(1) is equal to
$$\sum_i(-1)^i\sum_{k=0}^i
N{\frak q}^{n(\sigma)k}\sum_q(-1)^q
\text{Tr}(\sigma\circ \Gamma|
H^q(Y^{(i)},\Bbb Q_\ell)).$$

For an element $\sigma\in W^+_{\frak q}$
in the Weil group with
$n(\sigma)\ge 0$,
we define an endomorphism
$\sigma_{\text{geom}}$ of $Y^{(i)}$ to be
$\sigma_{\text{geom}}
=\sigma\circ (\text{abs. Frob.}
)^{[\Bbb F_{\frak q}:\Bbb F_p]\cdot n(\sigma)}$
for each $i$.
It is a geometric endomorphism of
a scheme $Y^{(i)}$ over the base field 
$\bar \Bbb F_{\frak q}$.
Since the absolute Frobenius acts trivially on
etale cohomology $H^q(Y^{(i)},\Bbb Q_\ell)$,
we have $\sigma_*=\sigma_{\text{geom} *}$
as an operator acting it.
Let $\Gamma_\sigma$ denote
the composite of $\sigma_{\text{geom}}$
with $\Gamma$ as an algebraic correspondence
and let $(\Gamma_\sigma,\Delta)$
be the intersection number.
We apply the Lefshetz trace formula,
to a proper smooth scheme $Y^{(i)}$
and an algebraic correspondence
$\Gamma_\sigma$.
Then we obtain
$$\sum_q(-1)^q
\text{Tr}(\sigma\circ \Gamma|
H^q(Y^{(i)},\Bbb Q_\ell))
=(\Gamma_\sigma,\Delta).$$

Next we compute the $p$-adic case.
For log cristalline cohomology,
we also have
the weight spectral sequence [M]
$$E_1^{i,j}=
\bigoplus_{k\ge \max(0,-i)}
H^{j-2k}_{\text{cris}}(Y^{(i+2k)}/W)(-k)\Rightarrow
H^{i+j}_{\text{log cris}}(Y/W).$$
Here the Tate twist $(-k)$ means that 
we replace the Frobenius $\varphi$ by
$p^k\varphi$.
Since the maps involved in
the definitions of the Hecke operators
are finite etale,
by the same argument as above,
we see that
the spectral sequence is compatible with
the action of Hecke operators,
endomorphisms and permutations.
It is also compatible
with the semi-linear action of 
the Galois group and the Frobenius operator.
Hence by modifying it in the same way
on the both side, it is also compatible
with the linear action of the Weil group.
For $\sigma\in W^+$ the modified action
$\sigma_*\circ F^{n(\sigma)[\Bbb F_{\frak q}:\Bbb F_p]}$
is the same as the action of
the geometric endomorphism
$\sigma_{\text{geom}}
=\sigma\circ (\text{abs. Frob.}
)^{[\Bbb F_{\frak q}:\Bbb F_p]\cdot n(\sigma)}$.
Hence the right hand side of (1) is equal to
$$\sum_i(-1)^i\sum_{k=0}^i
N{\frak q}^{n(\sigma)k}
\sum_q(-1)^q
\text{Tr}(\sigma_{\text{geom} *}\circ \Gamma|
H^q_{\text{cris}}(Y^{(i)}/W)).$$
Again by the Lefschetz trace formula [GM][Gr],
we have
$$\sum_q(-1)^q
\text{Tr}(\sigma_{\text{geom} *}\circ \Gamma|
H^q_{\text{cris}}(Y^{(i)}/W))
=(\Gamma_\sigma,\Delta).$$
Thus the both sides give the same answer
and the equality in (1) is proved.

Finally we prove the assertion (2),
the monodromy weight conjecture.
The algebraic correspondence $e$ in Lemma 3
acts as an projector on the
spectral sequences.
We consider the $e$-part of them.
We compute the $E_1$-terms of
the $e$-part.
We have $Y^{(k)}=\emptyset$
for $k>1$ since
the semi-stable model
$X_O$ is proper smooth
over the semi-stable model
$M_O$ of a curve.
Let $C$ denote the closed
fiber of the semi-stable model
$M_O$.
Then the disjoint union 
$C^{(0)}$ of the 
components is the same as the
normalization of $C$ and
the disjoint union 
$C^{(1)}$ of the 
their intersections is
the singular locus of $C$.
To describe the $E_1$-terms
we introduce some sheaves on
$C^{(i)}$.
For a place
$\lambda|\ell$ of $L$,
we define a smooth $L_\lambda$-sheaf
$\Cal F_\lambda^{(k)}$
to be
$\bigotimes_i
(\text{Sym}^{k_i-2}\otimes{\det}^{\frac{w_-k_i}2}
(e_iR^1a_*L_\lambda))\otimes
(e_0 R^1b_*L\lambda))^{\otimes (w-2)(g-1)}.$
It is the restriction of
the extension of
$\Cal F_\lambda^{(k)}$
on $M$ to $M_O$.
Similarly for a place
$\mu|p$ of $L$,
we define an $F$-isocrystal
$\Cal E_\lambda$.
We consider
$F$-isocrystals 
$R^1a_*O_{\text{crys}}\otimes_{W}\widehat{L_\mu^{nr}}$
and
$R^1b_*O_{\text{crys}}\otimes_{W}\widehat{L_\mu^{nr}}$
where the tensor product is taken 
as remarked after Claim 4.
We regard them
as an $O_D\otimes_{\Bbb Z}L$-module
and an $O_{E_0}\otimes_{\Bbb Z}L$-module
respectively.
Then we define
$\Cal E_\lambda$
to be
$$\bigotimes_i
(\text{Sym}^{k_i-2}\otimes{\det}^{\frac{w_-k_i}2}
(e_iR^1a_*O_{\text{crys}}\otimes_{W}\widehat{L_\mu^{nr}})\otimes
(e_0 R^1b_*O_{\text{crys}}\otimes_{W}\widehat{L_\mu^{nr}}
)^{\otimes (w-2)(g-1)}.$$
Then similarly as in Lemma 3,
we have
$e R^qc_*L_\lambda=\Cal F_\lambda^{(k)}$
if $q=q_0=(w-2)(2g-1)$ and $=0$ if otherwise.
Also we have
$e R^qc_*O_{\text{crys}}\otimes_{W}\widehat{L_\mu^{nr}}
=\Cal E_\mu$
if $q=q_0$ and $=0$ if otherwise.
By the same argument as in Lemma 3 using
the Leray spectral sequence,
we see that there are only 5 non-vanishing 
$E_1$-terms
$$\matrix
E_1^{-1,q_0+2}\quad& E_1^{0,q_0+2}&\\
                   & E_1^{0,q_0+1}&\\
                   & E_1^{0,q_0}&\quad E_1^{1,q_0}
\endmatrix$$
where $q_0=(2g-1)(w-2)$.
Each term is described as follows.
In the $\ell$-adic setting,
we have
$$E_1^{0,q_0+q}=H^q(C^{(0)},\Cal F_\lambda^{(k)}),\quad
E_1^{1,q_0}
=E_1^{-1,q_0+2}(1)=
H^q(C^{(1)},\Cal F^{(k)}).
$$
In the crystalline setting,
we replace 
$\Cal F_\lambda^{(k)}$
by $\Cal E_\lambda^{(k)}$.
The map $d_1^{-1,q_0+2}$
is the Gysin map and
$d_1^{1,q_0}$ is the restriction map.
By the Weil conjecture,
the eigenvalues of
a lifting of geometric Frobenius
acting on each $E_1$-term
$E_1^{i,j}$ are algebraic
integer purely of weight $j$
and the spectral sequence degenerates
at $E_2$-terms.
The monodromy operater $N$ on the limit
is induced by the canonical isomorphism
$N:E_1^{-1,q_0+2}(1)\to
E_1^{1,q_0}$ [RZ],[M].
Therefore the weight-monodromy conjecture
is equivalent to that
the isomorphism $N$ on the 
$E_1$-term induces an isomorphism on $E_2$-terms.
Thus it is reduced to show 

\proclaim{Claim 5}
Let $q_0=(2g-1)(w-2)$.
The canonical map
$$
N:\text{\rm Ker}(
E_1^{-1,q_0+2}\to E_1^{0,q_0+2})(1)
\to\text{\rm Coker}(E_1^{0,q_0} \to
\quad E_1^{1,q_0})$$
is an isomorphism.
\endproclaim

First we prove it in the case
where the multiweight $k$
is of the form $k=(2,2,\cdots,2,w)$.
In this case, the sheaves
$\Cal F_\lambda^{(k)}$ and $\Cal E_\mu^{(k)}$ are constant.
Let $I$ be the set of irreducible components
and $J$ be the set of singular points.
Then it is enough to show that
$\text{Ker}(\Bbb Q^J\to \Bbb Q^I)\to
\text{Coker}(\Bbb Q^I\to \Bbb Q^J)$
is an isomorphism.
It is proved easily by extending the scalar to
$\Bbb R$.

We assume
the multiweight $k$
is not of the form $k=(2,2,\cdots,2,w)$.
To show Claim,
we prove Proposition 1 below
in the next section.
To state it,
we introduce a terminology.
Take a sufficiently
small open compact subgroup
$K''$ such that $M''_{K''}$
has a proper smooth model $M''_{K'',O}$
and that $KH\subset K''$.
We say an component $C_i$ in $C$
is ordinary, if
it dominates a component
of the closed fiber of
$M''_{K'',O}$.
Otherwise, we say it is
supersingular.

\proclaim{Proposition 1}
Let $C_i$
be an ordinary
irreducible component of $C$.
Then we have
$$\align&
H^0(C_i,\Cal F^{(k)})=
H^2(C_i,\Cal F^{(k)})=0,\\
&H^0(C_i,\Cal E^{(k)})=
H^2(C_i,\Cal E^{(k)})=0
\endalign$$
unless $k=(2,2,\cdots,2,w)$.
\endproclaim

Proof will be given in the next section.

We show Claim admitting Proposition 1.
Let $\Sigma\subset M''_{K'',O}$ be the union
of the image of supersingular
components and of singular points.
Then for each $s\in \Sigma$
the sheaves
$\Cal F_\lambda^{(k)}$ and $\Cal E_\mu^{(k)}$ are constant
in the inverse image.
Let $I_x$ be the set of irreducible components
and $J_x$ be the set of singular points
in the inverse image.
Then it is reduced to that
$\text{Ker}(\Bbb Q^{J_x}\to \Bbb Q^{I_x})\to
\text{Coker}(\Bbb Q^{I_x}\to \Bbb Q^{J_x})$
is an isomorphism
which is proved in the same way as above.

\medskip

\subheading{9. Vanishing of $H^0$}

We prove Proposition 1.
First we restate it
in terms of the closed
fiber of the proper smooth model
$M'_{K',O_V}$ of $M'_{K'}$
and Tate-modules.
Let $K'\subset G'(\Bbb A_f)$
be a sufficiently small
open subgroup 
satisfying the condition 
(c) in the proof of
Lemma 4 (2) in Section 7:
$K'=\Bbb Z_p^\times\times GL_2(O_{F,\frak p})
\times K_p^{\prime {\frak p}}\times K^{\prime p}$.
Then as is recalled there,
Carayol has shown that 
$M'_{K'}$ has good reduction
and the abelian variety $A'_{K',\hat T}$
is extended to the proper smooth model
$M'_{K',O_{E,\frak q}}$.
Let $C$ be an irreducible component
of the geometric closed fiber
$M'_{K',O_{E,\frak q}}\otimes \bar \Bbb F_{\frak q}$.
We will define a smooth $\ell$-adic
sheaf $\Cal F^{*(k)}_\lambda$
and an $F$-isocrystal
$\Cal E^{*(k)}_\lambda$ on $C$
in a similar way as
$\Cal F^{\prime(k)}_\lambda$.
For a place
$\lambda|\ell$ of $L$,
we define a smooth $L_\lambda$-sheaf
$\Cal F^{*(k)}_\lambda$
to be
$$\bigotimes_i
\left(\text{Sym}^{k_i-2}(e_i T_\ell(A)\otimes_{\Bbb Z_\ell} L_\lambda)
\otimes(\det
(e_i T_\ell(A)\otimes_{\Bbb Z_\ell} L_\lambda)
)^{\otimes \frac{w_-k_i}2}\right).$$
Here the idempotents $e_i\in 
\text{End}_{M'}(A)\otimes L$
act on
$T_\ell(A)\otimes_{\Bbb Z_\ell} L_\lambda$
by the covariant functoriality
of Tate modules.
We define an $F$ crystal.
Let $\Cal T_p(A)$ denote
the $F$-crystal associated to
the $p$-divisible group
$A[p^\infty]$ on $M'$.
Let $\mu|p$ be a place of $L$.
We regard
the crystal
$\Cal T_p(A)\otimes_{W}\widehat{L_\mu^{nr}}$
as an $O_D\otimes_{\Bbb Z}L$-module
by the covariant functoriality as above.
For each $i$,
we define an $F$-isocrystal
$\Cal E_i$ to be
$e_i(\Cal T_p(A)\otimes_{W}\widehat{L_\mu^{nr}})$
and put 
$\Cal E^{*(k)}_\lambda
=\bigotimes_i
\left(\text{Sym}^{k_i-2}\Cal E_i
\otimes({\det}\Cal E_i)^{\otimes\frac{w_-k_i}2}\right)$.

\proclaim{Proposition 1'}
Let $C'\to C$ be a finite covering
of proper smooth curves
and assume the multiweight $k$
is not of the form $(2,2,\cdots,2,w)$.
For $\lambda|\ell\neq p$,
the pull-back to $C'$ of the smooth sheaf
$\Cal F^{*(k)}_\lambda$ has no
non-trivial (geometrically) 
constant subsheaf or quotient
smooth sheaf.
For $\mu|p$,
the pull-back to $C'$ of 
the underlying isocrystal
$\Cal E^{*(k)}_\lambda$ has no
non-trivial 
constant subisocrystal 
or quotient isocrystal.
\endproclaim

We show that Proposition 1' implies
Proposition 1.
Let $C_i$ be as in Proposition 1.
By the construction of
$\Cal F^{(k)}_\lambda$
and 
$\Cal E^{(k)}_\lambda$
on $C_i$, Proposition 1' implies
a similar statement where
we replace
$C', \Cal F^{*(k)}_\lambda$
and 
$\Cal E^{*(k)}_\lambda$ by
$C_i, \Cal F^{\prime(k)}_\lambda$
and 
$\Cal E^{\prime(k)}_\lambda$.
By definition of $H^0$ and
by Poincar\'e duality,
it implies the assertion in Proposition 1.

\demo{Proof of Proposition 1' for $\lambda\nmid p$}
First, we prove the $\ell$-adic case.
The argument is
similar to the 
proof of vanishing of $H^0$ and $H^2$
in the reduction of
the equality (1) in Claim 1
to that in Claim 2 given in Section 3.
It is enough to show that the
image of the action of $\pi_1(C)$ is
sufficiently large.
We show that the action on the Tate module
defines a surjection
$\pi_1(C)\to  SK'_\ell=
\text{Ker}(\nu:
K'_\ell \to \Bbb Z_\ell^\times
\times O_{E,\ell}^\times).$
Let $V$ denote the maximal unramified
extension of $E_{\frak q}$
and $M^{\prime+}_{K',O_V}$
be the connected component of 
the proper smooth model
whose closed fiber is $C$.
Since $T_\ell(A)$
is locally constant on 
$M^{\prime+}_{K',O_V}$
the map
$\pi_1(M^{\prime+}_{K',\bar V})\to 
(\hat O_D^\times)^p$
factors through a surjection
$\pi_1(M^{\prime+}_{K',\bar V})\to
\pi_1(M^{\prime+}_{K',O_V})\simeq
\pi_1(C)$.
Since
$\pi_1(M^{\prime+}_{K',\bar V})\simeq
\pi_1(M^{\prime+}_{K',\Bbb C})\simeq
SK'_\ell$,
we obtain the surjection.
The rest of the argument is
identical to that in loc.cit and 
we will not repeat it here.
\enddemo

To proceed to the crystalline case,
we recall the modular interpretation [C1]
due to Carayol of
the integral model of $M'$
over the integer ring $O=O_{E_{\frak p}}$.
Let $K'\subset G'(\Bbb A_f)$
be a sufficiently small
open subgroup 
satisfying the condition 
(c) in the proof of
Lemma 4 (2) in Section 7:
$K'=\Bbb Z_p^\times\times GL_2(O_{F,\frak p})
\times K_p^{\prime {\frak p}}\times K^{\prime p}$.
We take an order $O_D\subset D$ 
such that $K'\subset \hat O_D^\times$.
We take  an $\hat O_D$-lattice
$\hat T\subset D\otimes \Bbb A_f$.
We assume they satisfy the following conditions:

$O_D$ is stable under the involution $*$,

$O_D\otimes_{\Bbb Z} \Bbb Z_p$ is maximal in 
$D\otimes_{\Bbb Q}\Bbb Q_p$,

$\text{Tr} \psi(\hat T,\hat T)\subset \hat\Bbb Z$

and 
$\text{Tr} \psi(\hat T\otimes_{\hat O_E}O_{E_p},
\hat T\otimes_{\hat O_E} O_{E_p})
\to \Bbb Z_p$ is perfect.

\noindent We put
$\hat T^p=
\hat T\otimes_{\hat \Bbb Z}\prod_{q\neq p}\Bbb Z_q$.
It is a free $\hat O_E^p$-module of rank 4 
and has a symmetric bilinear from
$\text{Tr }\psi:\hat T^p\times \hat T^p\to \Bbb Z^p$.
We define a free $O_{F,p}^{\frak p}=
\prod_{\frak p'|p\neq \frak p}O_{F,\frak p'}$-module
$T_p^{\frak p}$ of rank 4 as follows.
By the isomorphism
$O_E\otimes\Bbb Z_p=\prod_{\frak p'|p}
(O_{F,\frak p'}\times O_{F,\frak p'}$,
we have direct sum decomposition
$\hat T\otimes \Bbb Z_p=\prod_{\frak p'|p}
(\hat T_{\frak p'_1}\times T_{\frak p'_2})$.
Here the first factors correspond to the
embedding $O_{E_0}\to \Bbb Z_p$ fixed in Section 7.
We put $T_p^{\frak p}=
\prod_{\frak p'|p,\neq \frak p}
\hat T_{\frak p'_1}$.

For an $O_D$-abelian scheme $A$ on
an $O_{E,\frak q}$-scheme $S$,
we define direct summands $\text{Lie}^2A$
and $\text{Lie}^{1,2}A$ of $\text{Lie}A$
similarly as in Section 5.
We define $T_p^{\frak p}(A)$ similarly as
$\hat T_p^{\frak p}$ as above.
On the category of schemes over $O=O_{E_{\frak p}}$,
there is a proper smooth model
$M'_{K',O}$ over $O$ of $M'_{K'}$ representing the functor
$S\mapsto \{\text{isomorphism class of }(A,\theta,\bar k)\}$
where 
\roster\item
$A$ is an $O_D$-abelian schemes of dimension $4g$
such that $\text{Lie}^2 A=\text{Lie}^{1,2} A$
and it is a locally free $O_S$-module of rank 2.
\item
$\theta\in \text{Hom}(A,A^*)^{sym}$ 
is an $O_D$-polarization of $A$.
\item
$\bar k=\bar k^p_{\frak p}\times \bar k^p$ is a 
pair of a $K_p^{\prime \frak p}$-equivalent class of
a $\prod_{{\frak p}'|p,\neq {\frak p}}
O_{D_{\frak p'}}$-isomorphism 
$k^{\frak p}_p:T_p^{\frak p}(A)\to T_p^{\frak p}$
and a $K^{\prime p}$-equivalent class of
a $\prod_{{\frak p}'\nmid p}O_{D_{\frak p}'}
$-isomorphism $k^p:T^p(A)\to T^p$
such that there exists a $ \hat\Bbb Z^p=
\prod_{q\neq p}\Bbb Z_q$-isomorphism $k'$
making the diagram 
$$\matrix
T^p(A)\times T^p(A)&@>{(1,\theta_*)}>> 
T^p(A)\times T^p(A^*)@>>> &\hat\Bbb Z^p(1)\\
@V{k\times k}VV @VV{k'}V\\
T^p\times T^p &@>>{Tr \psi}>& \hat\Bbb Z^p
\endmatrix$$
commutative.
\endroster
In (3), the $O_E\otimes \hat\Bbb Z^p$-module $T^p(A)$ is free of rank 4
and, by the condition (1),
the $\prod_{{\frak p}'|p,{\frak p}'\neq {\frak p}}
O_{E_{\frak p}'}$-module $T_p^{\frak p}(A)$ is also free of rank 4.
As is shown in [C1],
the generic fiber $M'_{K,O}\otimes_OE_{\frak p}$
represents the restriction of the functor
$M'_K$ to the schemes over $E_{\frak p}$.
Hence the smooth proper scheme
$M'_{K,O}$ is a model of
the base change $M'_K\otimes_EE_{\frak p}$
and the universal abelian scheme $A$
is a unique extension on $M'_{K,O}$
of the pull-back.

We state Lemma 5 and 6 on the $p$-divisible group
$A[p^\infty]$ on $C$.
We will deduce Proposition 1' in crystal case
from the Lemmas.
As in [C1] 2.6.3,
we put $$T_p^{\frak p}(A)=
T_p(A)\otimes_{O_{E,p}}
\prod_{\frak q'|p,\nmid \frak q_{0},\nmid \frak p}
O_{E,\frak q'}.$$
We identify
$\prod_{\frak q'|p,\nmid \frak q_{0},\nmid \frak p}
O_{E,\frak q'}=
\prod_{\frak p'|p,\neq \frak p}
O_{F,\frak p}=O_{F,p}^{\frak p}$
and regard
$T_p^{\frak p}(A)$
as an 
$O_{B,p}^{\frak p}=
\prod_{\frak p'|p,\neq \frak p}
O_{B,\frak p}$-module.
By the modular interpretation recalled above,
it is a smooth etale sheaf on 
the proper scheme $M'_{K',O}$ of
$O_{B,p}^{\frak p}=
\prod_{\frak p'|p,\neq \frak p}
O_{B,\frak p}$-modules
of rank 1.
Let $\frak q_2|\frak p,\neq \frak q$ be 
the other prime ideal of $O_E$
dividing $\frak p$.
We identify
$O_{D,\frak q_2}=
O_{B,\frak p}$ and
take an isomorphism
$O_{B,\frak p}\simeq
 M_2(O_{F,\frak p})$.
Let $e\in O_{D,\frak q_2}$ be the idempotent
corresponding to
$\pmatrix 1& 0\\ 0&0
\endpmatrix
\simeq M_2(O_{F,\frak p})$.
Similarly as in [C1] 5.4,
let $\Bbb E_\infty$
be the $p$-divisible group
$$\Bbb E_\infty
=e(A[p^\infty]
\otimes_{O_E\otimes \Bbb Z_p}
O_{E,\frak q_2}).$$
In the terminology of
[C1] Appendix 1,
it is an $O_{F,\frak p}$-divisible group
of height 2.
As a $p$-divisible group,
it is of height 
$2[F_{\frak p}:\Bbb Q_p]$
and of dimension $1$.
Let $\Cal E_0$ and
$\Cal T^{\frak p}$
be the $F$-crystals
associated to
the $p$-divisible group
$\Bbb E_\infty$
and to the Tate module
$T_p^{\frak p}(A)$
respectively.

The $F$-isocrystal
$\Cal E^{*(k)}_\mu
=\bigotimes_i
(\text{Sym}^{k_i-2}\Cal E_i
\otimes({\det}\Cal E_i)^{\frac{w_-k_i}2})$
is related to
them in the following way.
We regard $\Cal E_0$ and
$\Cal T^{\frak p}$
as an $O_{F,\frak p}$-module
and an $O_{F,p}^{\frak p}$-module
respectively by the covariant functoriality.
Let $I_1\subset I=\{\tau_i:F\to L\}$ be 
the subset $I_1=\{\tau_i:F_{\frak p}\to L_{\mu}\}$.
Then for $i\in I_1$,
the $F$-isocrystal $\Cal E_i$
is isomorphic to $\Cal E_0\otimes_{O_{F,\frak p}}
\widehat{L_\mu^{nr}}$
with respect to
$\tau_i:O_{F,\frak p}\to \widehat{L_\mu^{nr}}$.
Here we identify 
$O_{F,\frak p}$ with
$O_{E,\frak q}$.
For $i\in I-I_1$,
we take an isomorphism
$B\otimes_FL\simeq M_2(L)$
for the tensor product with respect to $\tau_i$
and let $e$ be the idempotent
corresponding to 
$\pmatrix 1&0\\0&0\endpmatrix.$
Then the $F$-isocrystal $\Cal E_i$
is isomorphic to $e
(\Cal T^{\frak p} \otimes_{O_{F,p}^{\frak p}}
\widehat{L_\mu^{nr}})$
with respect to
$\tau_i:O_{F,p}^{\frak p}\to \widehat{L_\mu^{nr}}$.
Here we also identify
$O_{E,\frak q'}=O_{F,\frak p'}$
for primes $\frak p'|p,\neq \frak p$
and $\frak q'|\frak p',\nmid \frak q_0$.

It is snown in [Ca] (6.7), (9.4.3) that
there exists a finite nonempty set $\Sigma\subset C$
of closed points satisfying the following condition:

At each point in $\Sigma$,
the $p$-divisible group $\Bbb E_\infty$ is connected.
On the complement
$U=C-\Sigma$,
the $p$-divisible group $\Bbb E_\infty$
is an extension of an etale $p$-divisible group
$\Bbb E_\infty^{\text{et}}$
by an connected $p$-divisible group
$\Bbb E_\infty^\circ$.

\noindent
We call a point in $\Sigma$ a supersingular point
and a point in $U$ an ordinary point.
The $p$-divisible groups
$\Bbb E_\infty^{\text{et}}$
and 
$\Bbb E_\infty^\circ$
have natural structures of $O_{F,\frak p}$-modules.
The Tate module
$T(\Bbb E_\infty^{\text{et}})$
is a smooth sheaf of $O_{F,\frak p}$-modules of rank 1.

\proclaim{Lemma 5} 
The morphism 
$\pi_1(U)\to O_{F, \frak p}^\times\times
O_{B,p}^{\frak p\times}$ defined
by the smooth sheaf 
$T_p(\Bbb E_\infty^{\text{et}})\times
T_p^{\frak p}(A)$ of $O_{F, \frak p}\times
O_{B,p}^{\frak p}$-modules of rank 1
defines a surjection
$$\pi_1(U)\to O_{\frak q}^\times\times
SK^{\prime {\frak q}}.$$
\endproclaim

Let $\Cal E'_0$ and $\Cal E''_0$
be the $F$-crystal associated to
the $p$-divisibles 
$\Bbb E_\infty^{\text{et}}$ and
$\Bbb E_\infty^\circ$ on the ordinary locus
$U$ respectively.
The restriction of $\Cal E_0$ on $U$
is an extension
$$0@>>> \Cal E'_0 @>>> \Cal E_0 @>>> \Cal E''_0 @>>> 0.$$
as $\Bbb E_\infty$ is an extension.

\proclaim{Lemma 6}
The extension of the underlying isocrystal
$$0@>>> \Cal E'_0 \otimes\Bbb Q_p
@>>> \Cal E_0 \otimes\Bbb Q_p
@>>> \Cal E''_0 \otimes\Bbb Q_p@>>> 0$$
is non-trivial.
\endproclaim

\demo{Proof of Proposition 1' for $\lambda\nmid p$}
The argument is similar to that in [Cr].
First we prove it admitting Lemmas 5 and 6.
It is sufficient to show that,
on the inverse image $U'\subset C'$
of the ordinary locus $U$,
the restriction of $\Cal E^{*(k)}$ 
has no non-constant subisocrystal
or quotient isocrystal.
Before starting proof,
note that an $F$-(iso)crystal is
constant if and only if
the underlying (iso)crystal is constant.
In fact, if the underlying (iso)crystal $\Cal E$ is constant,
the Frobenius pull-back $F^*\Cal E$
and the Frobenius map $F:F^*\Cal E\to \Cal E$
defining the structure of
$F$-(iso)crystal is constant.
The only if part is trivial.

We put
$r=[F_{\frak q}:\Bbb Q_p]$
and $I_1=\text{Hom}(F_{\frak q},L_\mu)=
\{\tau_1,\ldots,\tau_r\}
\subset I=\text{Hom}(F,L)=
\{\tau_1,\ldots,\tau_g\}.$
We define a decreasing filtration on
the restriction of $\Cal E^{*(k)}$ on $U$
with multiindex $\Bbb Z^{I_1}$ as follows.
On $\Cal E_0$, 
we define a filtration $F^\bullet $ on $\Cal E_0$ by
$F^0\Cal E_0=\Cal E_0, F^1\Cal E_0=\Cal E'_0,
F^2\Cal E_0=0$.
For each $i\in I_1$,
it induces a filtration on
$\Cal E_i$ and hence on $\text{Sym}^{k_i-2}\Cal E_i$
by the isomorphism
$\Cal E_i\simeq
\Cal E_0\otimes_{O_{F,\frak p}}L_\mu$.
Taking symmetric powers and tensor product,
we obtain a filtration on
$\Cal E^{*(k)}=
\bigotimes_{i\in I}
(\text{Sym}^{k_i-2}\Cal E_i\otimes
(\det \Cal E_i)^{\otimes \frac{w-k_i}2})$.
We consider the graded piece
$Gr_F^q \Cal E^{*(k)}
=F^q \Cal E^{*(k)}/
\sum_{q'>q}F^{q'} \Cal E^{*(k)}$
for each $q=(q_1,\ldots,q_r)\in \Bbb Z^{I_1}$.

We deduce from Lemma 5 that
the isocrystal $Gr_F^q\Cal E^{(k)}$
has no constant sub\-iso\-crys\-tal
or quotient isocrystal
except for at most one multi-index
$q=(q_1,\ldots, q_r)$ satisfying
$(k_1,\cdots, k_g)
=(2q_1+2,\ldots, 2q_r+2,2,\ldots,2)$.
In the exceptional case,
we will see that
the graded piece is in fact constant.
We compute the graded pieces.
The graded pieces are computed as
$$\align
&Gr_F^q \Cal E^{(k)}\\
=&\bigotimes_{i\in I_1}
((\det \Cal E_i)^{\otimes \frac{w-k_i}2+q_i}\otimes 
(Gr^0_F\Cal E_i)^{\otimes k_i-2-2q_i})\otimes 
\bigotimes_{i\in I-I_1}
(\text{Sym}^{k_i-2} \Cal E_i
\otimes(\det \Cal E_i)^{\otimes \frac{w-k_i}2})
\endalign
$$
for $0\le q_i\le k_i-2$ for $i\in I_1$
and is 0 otherwise.
By the Weil-pairing of
Drinfeld basis [C1] 9.2,
the determinant isocrystal
$\det \Cal E_i$
is geometrically constant
for $i\in I_1$.
Similarly but more easily,
$\det \Cal E_i$ is also constant for $i\in I-I_1$.
Therefore it is sufficient to show that
the isocrystal
$\bigotimes_{i\in I_1}
(Gr^0_F\Cal E_i)^{\otimes k_i-2-2q_i}
\otimes \bigotimes_{i\in I-I_1}
\text{Sym}^{k_i-2} \Cal E_i$
has no non-trivial
constant subisocrystal or quotient isocrystal
unless $(k_1,\cdots, k_g)
=(2q_1+2,\cdots,2q_r+2,2,\cdots,2)$.

The $F$-isocrystals
$Gr^0_F\Cal E_i$ for $i\in I_1$ and
$\Cal E_i$ for $i\in I -I_1$ are
defined by smooth $p$-adic etale sheaves
on $U$.
Let $\Cal L_i$ and $\Cal F_i$ be
the corresponding smooth $p$-adic sheaves.
Since $F$-isocrystal is constant if and only
if the underlying crystal is constant,
it is reduced to show that
the smooth $p$-adic sheaf
$\bigotimes_{i\in I_1}
\Cal L_i^{\otimes k_i-2-2q_i}
\otimes \bigotimes_{i\in I-I_1}
\text{Sym}^{k_i-2} \Cal F_i$
is an irreducible
smooth $p$-adic etale sheaf.
It follows from the surjectivity, Lemma 5,
of the map
$\pi_1(U)\to SK'_p$
by the same argument as in the $\ell$-adic
case.

We complete the proof by using Lemma 6.
We assume that there exists a non-trivial
constant subisocrystal of $\Cal E^{*(k)}$
for $(k_1,\cdots, k_g)
\neq (2,\cdots,2)$ and deduce a contradiction.
The proof for the quotient is similar
and is omitted.
By the study of the graded pieces above,
the proof is completed
except for the case where
$k_i$ are even for $i\in I_1$
and $k_i=2$ for $i\in I-I_1$.
We put
$(k_1,\cdots, k_g)
=(2q_1+2,\ldots, 2q_r+2,2,\ldots,2)$
and assume $q=(q_1,\ldots, q_r)\neq 0$.
By the computation of the graded pieces,
if we had non-trivial
constant subisocrystal,
it should be contained in
$F^q \Cal E^{*(k)}$ and 
mapped isomorphically to
$Gr^q \Cal E^{*(k)}$.
Namely, the extension
$F^q \Cal E^{*(k)}$
of 
$Gr^q \Cal E^{*(k)}$ is split.
Take an index $i\in I_1$ such that
$q_i>0$
and let $q',q''$ be the multi-index
obtained from $q$ by replacing
$q_i$ be $q_i+1,q_i+2$
respectively.
Then the extension
$$0@>>>
Gr^{q'} \Cal E^{*(k)}
@>>>
F^q \Cal E^{*(k)}/
F^{q''} \Cal E^{*(k)}
@>>>
Gr^q \Cal E^{*(k)}@>>>0$$
is also split.
Its extension class
is $q_i$-times the class 
of extension in Lemma 6 and hence
is non-zero.
Thus we get a contradiction.
We have proved that Lemmas 5 and 6
implies Proposition 1'.
\enddemo

We prove Lemmas 5 and 6
to complete the proof of Proposition 1'
hence of Theorems 1 and 2.
We prove Lemma 5 using a supersingular point
which exists by [Ca] (9.4.3).
Lemma 6 will be proved using an ordinary point.

\demo{Proof of Lemma 5}
Since $T_p^{\frak p}(A)$ is
smooth on the proper smooth model
$M'_{K',O_{E,\frak q}}$,
the same argument as in the $\ell$-adic case
shows that we have a 
surjection $\pi_1(C)\to SK_p^{\prime\frak p}$.
Take a supersingular point $x\in \Sigma\neq\emptyset$
and let $I_x$ denote the inertia group.
It is enough to show that
the restriction $I_x\to O_{\frak q}^\times$
is surjective.
Let $U_n$ be the finite etale covering
$U_n=\text{Isom}
(O_{F,\frak p}/\frak p^n,\Bbb E_n^{et})$
of $U$ trivializing
the ${\frak p}^n$-torsion part $\Bbb E_n^{et}$.
Here an isomorphism means
an isomorphism of 
$O_{F,\frak p}/\frak p^n$-group schemes.
The covering $U_n$
is an analogue of an Igusa curve.
It is sufficient to show that
$U_n$
is totally ramified at a supersingular point.
Namely we show the following.

\proclaim{Lemma 7} 
Let $K_x$ denote the completion
of the function field of $C$ at 
a supersingular point $x$.
Then the base change 
$U_n\times_C\text{Spec }K_x$
is the spectrum of a totally ramified
extension of $K_x$.
\endproclaim

\demo{Proof}
Let $E$ denote the formal group associated
to the $p$-divisible group $\Bbb E_\infty$
over the completion
$\hat C=\text{Spec }\hat O_{C,x}$.
Let $\pi$ be a prime element of
$O_{F,\frak p}$.
For an integer $n$,
let $E^{(n)}$ denote the base change of
$E$
by the $N\frak p^n$-th power Frobenius
and $F^n:E\to E^{(n)}$ be the
$N\frak p^n$-th power relative Frobenius
over $\hat C$.
Then the multiplication $[\pi^n]:
E\to E$ is factorized as
$[\pi^n]=V^n\circ F^n$
for a map
$V^n:E^{(n)}\to E$.
Outside the closed point $x$,
the map $V^n$ is etale and hence
$\text{Ker}V^n$ is
a finite flat group scheme over $\hat C$
extending the etale quotient $\Bbb E_n^{et}$
on the generic point.
Let $C_n=(\text{Ker}V^n)^\times$
be the scheme of
$O_{F,\frak p}/\frak p^n$-basis of 
$\text{Ker} V^n$ in the sense of Drinfeld.
Namely,
it is a closed subscheme of
$\text{Ker} V^n$ representing the functor
$$R\mapsto\{s\in \text{Ker} V^n(R)|
\sum_{a\in O_{F,\frak p}/\frak p^n}
[as]=\text{Ker} V^n\text{ as a divisor in }
E^{(n)}_R\}$$
for a ring over $\hat O_{C,x}$.
Outside the closed point,
the scheme $C_n$ is the same as 
the base change of $U_n$.
Therefore,
it is sufficient to show that
$C_n$ is regular and
the inverse image of the closed point
$x$ by $C_n\to \hat C$
contains only one point.
The second assertion is clear
since $C_n$ is a closed subscheme of
a local scheme
$\text{Ker} V^n$.
We show that the intersection 
$C_n\cap [0]$ of $C_n$
with the zero-section $[0]$
of the formal group $E^{(n)}$
is equal to $\text{Spec }\kappa(x)$.
This implies that $C_n$ is regular
since the zero section is a divisor in 
$E^{(n)}$.

Let $R=\Gamma(C_n\cap [0], O)$.
It is an Artin $\hat O_{C,x}$-algebra.
It is sufficient to show that
the surjection
$\hat O_{C,x}\to R$
factors through the surjection
$\hat O_{C,x}\to \kappa(x)$.
By the assumption,
the zero-section is an 
$O_{F,\frak p}/\frak p^n$-basis of
$\text{Ker} V^n$.
Hence, we have
$\text{Ker} [\pi^n]=
\text{Ker} F^{2n}$
on $R$ and an isomorphism
$E_R\simeq E_R^{(2n)}\simeq
E_R^{(2mn)}$ for $m\ge 1$.
Since $R$ is Artinean,
for sufficiently large $m$,
the map $a\to a^{N\frak q^{2mn}}$
factors through $R\to \kappa(x)\to R$
and we obtain
$E_R\simeq E_R^{(2mn)}
\simeq E_x\otimes_{\kappa(x)}R$.
This means that 
$\hat O_{C,x}\to R$ factors through $\kappa(x)$
since $\Bbb E_\infty$ over $\hat O_{C,x}$ is
the universal deformation of 
$\Bbb E_\infty|_x$,
Proposition 5.4 [C1].
Thus we have proved Lemma 7 and hence Lemma 5.
\enddemo\enddemo

To prove Lemma 6,
we show

\proclaim{Lemma 8}
Let $\hat C=\text{Spec }\hat O_{C,x}$
be the completion at an ordinary
closed point $x\in U$.
Let 
$[\Bbb E]\in \text{\rm Ext}^1(\Bbb E^{et},\Bbb E^\circ)$
be the class of $\Bbb E$
as an extension of
$O_{F,\frak p}$-divisible groups
on $\hat C$.
Then the class
$[\Bbb E]$ is not torsion.
\endproclaim

We derive it from the following
statement proved in
[C1] Proposition 5.4,
App. Th\'eor\`eme 3.

\proclaim{Lemma 9} 
On the completion $\hat C$
at an ordinary closed point,
the connected part $\Bbb E^\circ$
is isomorphic to
the pull-back of the Lubin-Tate
formal group.
The etale part 
$\Bbb E^{et}$ is isomorphic to
the constant $O_{F,\frak p}$-divisible group
$F_{\frak p}/O_{F,{\frak p}}$.
The completion $\hat C$
pro-represents the functor
$R\mapsto 
\text{\rm Ext}_R(F_{\frak p}/O_{F,{\frak p}},
\Bbb E_0)=\Bbb E_0(R)$ 
on the category of Artin $\Bbb F_{\frak p}$-algebras
$R$ together with a surjection
$R\to \kappa(x)$.
It is isomorphic
to $\Bbb E^0$ as a formal scheme.
The extension $\Bbb E$ on
$\hat C=\Bbb E^0$ is identified with
the universal extension.
\endproclaim

\demo{Proof of Lemma 8}
We identify the formal schemes $\Bbb E^0=\hat C$.
By Lemma 9, the universal extension $\Bbb E$
corresponds to the identity $C\to \Bbb E^0$.
Hence it is the universal section of 
the formal group $\Bbb E^0$ and is not torsion.
\enddemo

\demo{Proof of Lemma 6}
It is enough to prove that
the restriction
to the completion at an ordinary closed point
is not the trivial extension.
Since the $p$-divisible groups
$\Bbb E^\circ$ and $\Bbb E^{et}$ are constant
on $\hat C$,
the $F$-isocrystals
$\Cal E'\otimes\Bbb Q_p, \Cal E''\otimes\Bbb Q_p $
and hence their underlying crystals
are constant there.
If the extension of the underlying
isocrystal was trivial,
the underlying isocrystal and hence the
$F$-isocrystal
$\Cal E\otimes\Bbb Q_p$
would be constant.
It means that
the extension class
$[\Bbb E]\in \text{\rm Ext}^1(\Bbb E^{et},\Bbb E^\circ)$
is torsion and contradicts with Lemma 8.
\enddemo

Thus the proof of Proposition 1'
and hence of Theorems 1 and 2
are now complete.

\enddocument